\documentclass[11pt]{amsart}
\usepackage{amsmath,amssymb,amsthm,mathrsfs}
\usepackage[margin=1in]{geometry}
\usepackage{tikz}
\usepackage{hyperref}

\newtheorem{theorem}{Theorem}[section]
\newtheorem{proposition}[theorem]{Proposition}
\newtheorem{lemma}[theorem]{Lemma}
\newtheorem{corollary}[theorem]{Corollary}
\newtheorem{conjecture}[theorem]{Conjecture}
\newtheorem{introtheorem}{Theorem}

\newtheorem{introconjecture}[introtheorem]{Conjecture}
\newtheorem{introcorollary}[introtheorem]{Corollary}
\theoremstyle{definition}
\newtheorem{definition}[theorem]{Definition}
\newtheorem{example}[theorem]{Example}
\newtheorem{problem}[theorem]{Problem}

\theoremstyle{remark}
\newtheorem{remark}[theorem]{Remark}

\usepackage{ulem}
\usepackage{color}

\newcommand{\bS}{\mathbb{S}}
\newcommand{\kk}{\Bbbk}
\newcommand{\Ord}{\mathrm{Ord}_+}
\newcommand{\Vect}{\mathrm{Vect}}
\newcommand{\Sym}{\mathrm{Sym}}

\newcommand{\lm}{\mathrm{lm}}
\newcommand{\lc}{\mathrm{lc}}
\newcommand{\st}{\mathrm{st}}
\newcommand{\sgn}{\mathrm{sgn}}
\newcommand{\cA}{\mathscr{A}}
\newcommand{\cB}{\mathcal{B}}
\newcommand{\Cl}{\mathrm{Cl}}

\newcommand{\cI}{\mathscr{I}}
\newcommand{\cJ}{\mathscr{J}}
\newcommand{\cM}{\mathscr{M}}
\newcommand{\AMR}{\cA_{MR}}
\newcommand{\cS}{\mathscr{S}}

\newcommand{\cW}{\mathcal{W}}
\newcommand{\Bij}{\mathrm{Bij}}
\newcommand{\Av}{\mathrm{Av}}
\newcommand{\id}{\mathrm{id}}

\newcommand{\lex}{\mathrm{lex}}
\newcommand{\Mat}{\mathrm{Mat}}
\newcommand{\Ess}{\operatorname{Ess}}

\newcommand{\End}{\operatorname{End}}
\newcommand{\ad}{\mathrm{ad}}
\newcommand{\dg}{\mathrm{d}}

\newcommand{\gridlines}[2]{\draw[gray!55,very thin] (0,0) grid (#1,-#2);}
\newcommand{\shaderows}[2][gray!22]{\foreach \x/\y in {#2}{\fill[#1] (0,1-\x) rectangle (\y,-\x);}}
\newcommand{\shadecell}[3][gray!50]{\fill[#1] (#3-1,-#2) rectangle (#3,1-#2);}
\newcommand{\pt}[2]{\fill (#2-0.5,0.5-#1) circle (0.16);}
\newcommand{\opt}[2]{\draw[semithick,fill=white] (#2-0.5,0.5-#1) circle (0.16);}
\newcommand{\rays}[3]{\draw[semithick] (#3,0.5-#1) -- (#2-0.5,0.5-#1) -- (#2-0.5,-#3);}
\newcommand{\crosscell}[2]{\draw[thick] (#2-0.8,0.2-#1) -- (#2-0.2,0.8-#1); \draw[thick] (#2-0.8,0.8-#1) -- (#2-0.2,0.2-#1);}
\newcommand{\wit}[2]{\draw[thick] (#2-0.5,0.5-#1) circle (0.34);}

\author{Anton Khoroshkin}
\address{Anton Khoroshkin: Department of Mathematics, University of Haifa, Mount Carmel, 3498838, Haifa, Israel}
\email{khoroshkin@gmail.com}

\title{Commutative bishuffle algebras and Pattern Avoidance}

\begin{document}

\begin{abstract}
The main purpose of this paper is to introduce the notion of a commutative bishuffle algebra, which linearizes classical pattern avoidance. The free commutative bishuffle algebra on one generator has a basis of permutations, divisibility of monomials is pattern containment, and the classes of permutations avoiding given patterns are the monomial bases of its quotients by monomial ideals. This free algebra is not Noetherian, but we show that its quotient by the antisymmetrizer $e_k$ is Noetherian for every $k$, and we conjecture that the same holds for its quotient by any nonzero ideal stable under the symmetric groups. We develop Gr\"obner bases for ideals of the free algebra and prove that the elements $e_k\oplus\tau$, with $\tau$ ranging over any set of permutations, form a Gr\"obner basis for the lexicographic ordering and for its opposite. This gives an algebraic form of the theorem of Backelin, West and Xin: the permutations avoiding all $k\cdots21\oplus\tau$ and those avoiding all $12\cdots k\oplus\tau$ are two monomial bases of one and the same quotient.
\end{abstract}

\maketitle

{\small
\tableofcontents
}
\setcounter{section}{-1}

\section{Introduction}

The main motivation for this paper is to advertise the notion of a \emph{commutative bishuffle algebra}, which linearizes the notion of classical pattern avoidance. We hope that the Noetherianity statements and questions and the Gr\"obner basis ideas discussed below can shed new light on well-known results on pattern avoidance. The author is not an expert in pattern avoidance, and this paper does not aim at striking new theorems in that area. Rather, we want to suggest that the extra algebraic structure may offer a new point of view on existing theorems.

Let us recall that a permutation $w\in\bS_n$ \emph{contains} a permutation $p\in\bS_k$ if there exist indices $1\le i_1<\dots<i_k\le n$ such that for all $1\le a,b\le k$ we have $w(i_a)<w(i_b)$ if and only if $p(a)<p(b)$. In other words, the permutation matrix of $p$ is obtained from that of $w$ by deleting some rows and the same number of columns. Otherwise $w$ \emph{avoids} $p$. For a set $P$ of permutations, $\Av_n(P)$ is the set of permutations of $[n]$ avoiding all elements of $P$. Two sets $P$ and $P'$ are \emph{Wilf-equivalent} if $\#\Av_n(P)=\#\Av_n(P')$ for all $n$.

If one thinks of {\it containment} as {\it divisibility} of monomials in a certain algebraic structure, one immediately arrives at the notion of what we call a \emph{commutative bishuffle algebra} (see Definitions~\ref{def::bishuffle} and~\ref{def:commutative} for details).
Roughly speaking, a sequence of vector spaces $\cA(n)$ with associative products
\[
\star_{\sigma,\tau}\colon\cA(m)\otimes\cA(n)\to\cA(m+n)
\]
indexed by pairs $(\sigma,\tau)$ of $(m,n)$-shuffles (that is, by pairs of decompositions of $[m+n]$ into a subset of size $m$ and a subset of size $n$) is called a \emph{bishuffle algebra}. It is \emph{commutative} if the product does not change when the two factors and the two parts of both decompositions are exchanged:
\[
u\star_{\sigma,\tau}v=v\star_{\sigma^{op},\tau^{op}}u.
\]
The free commutative bishuffle algebra $\cS$ on one generator of arity one is very concrete. Its component $\cS(n)=\kk\bS_n$ has a basis of permutation matrices of size $n$, and $u\star_{\sigma,\tau}v$ is obtained from the block-diagonal matrix with the blocks $u$ and $v$ by interleaving its rows according to $\sigma$ and, independently, its columns according to $\tau$ (Figure~\ref{fig:matrices}). Hence a permutation $w$ is divisible by $p$ in $\cS$ if and only if $w$ contains $p$. For a set $P$ of permutations, the quotient of $\cS$ by the monomial ideal $(P)$ has a basis consisting of the permutations avoiding $P$. Thus pattern avoidance problems become problems about one-generated commutative bishuffle algebras with monomial relations, just as consecutive pattern avoidance is the theory of monomial ideals of the free shuffle algebra \cite{DK}.

The advantage of the linear setting is that one may also consider ideals that are not monomial, and this gives information about monomial ideals as well.
For example, if an ideal $\cI$ of $\cS$ has two different monomial degenerations, then the corresponding two sets of leading monomials are Wilf-equivalent, because the permutations avoiding either of them form a basis of the same quotient $\cS/\cI$ (Proposition~\ref{prop:wilf}).

Another topic of this paper is the Noetherianity of finitely generated commutative bishuffle algebras.
We call a commutative bishuffle algebra \emph{Noetherian} if all its ideals are finitely generated. The algebra $\cS$ is not Noetherian: the containment order on permutations has infinite antichains \cite{SB,AMR,Vatter}, and such an antichain generates a monomial ideal that is not finitely generated (Proposition~\ref{prop:antichain}). The situation changes after passing to the quotient by the antisymmetrizer $e_k=\sum_{\sigma\in\bS_k}\sgn(\sigma)\,\sigma$.

\begin{introtheorem}[Theorem~\ref{thm:noeth-ek}]\label{introthm:noeth}
For every $k\ge1$ the commutative bishuffle algebra $\cS/(e_k)$ is Noetherian.
\end{introtheorem}

This cannot be seen at the level of monomials: for $k\ge3$ the permutations avoiding $k\cdots21$, which form a basis of $\cS/(e_k)$, contain an infinite antichain. The proof uses another commutative bishuffle algebra, $\cB_V=\bigoplus_n\End(V^{\otimes n})$ with $\dim V=k-1$. By Schur--Weyl duality, $\cS/(e_k)$ is the subalgebra of $GL(V)$-invariants of $\cB_V$. The monomials of $\cB_V$ are pairs of words in a finite alphabet, so they are well-quasi-ordered by Higman's lemma \cite{Higman}, and $\cB_V$ is Noetherian. The passage to the invariants is Hilbert's classical argument with the Reynolds operator, and the whole proof follows the strategy of Sam and Snowden \cite{SS}.

The algebra $\cS/(e_k)$ is an example of an \emph{$\bS$-symmetric} commutative bishuffle algebra, that is, of a quotient of $\cS$ by an ideal whose components are two-sided ideals of the group algebras $\kk\bS_n$. These algebras are classified by the upper sets of Young's lattice (Proposition~\ref{prop:symideals}), and they are the bishuffle analogues of twisted commutative algebras (Remark~\ref{rem:twisted}).

\begin{introconjecture}[Conjecture~\ref{conj:noeth-hooks}]\label{introconj:noeth}
Every $\bS$-symmetric commutative bishuffle algebra other than $\cS$ itself is Noetherian.
\end{introconjecture}

Every such algebra is a quotient of one of the algebras $\cA_{k,l}$. By the super Schur--Weyl duality of Sergeev and of Berele and Regev \cite{Sergeev,BR}, $\cA_{k,l}$ is the algebra of $\mathfrak{gl}(k|l)$-invariants in $\cB_V$ for the super vector space $V=\kk^{k|l}$. Theorem~\ref{introthm:noeth} is the case $l=0$. The first open case is $\cA_{1,1}$, since for $k,l\ge1$ the general linear supergroup is not reductive and there is no Reynolds operator. Finiteness questions of this kind also appear in the Specht problem on T-ideals \cite{Specht}, solved by Kemer \cite{Kemer87,Kemer91}.

In Section~\ref{sec:GB} we develop the theory of Gr\"obner bases for ideals of $\cS$, including the Buchberger criterion. For example, $e_k$ is a Gr\"obner basis of $(e_k)$ for every monomial ordering, with the leading monomial $k\cdots21$ or $12\cdots k$. This statement is equivalent to Schensted's enumeration of the permutations without long monotone subsequences \cite{Schensted} (Theorem~\ref{thm:en} and Remark~\ref{rem:rsk}).

As an application of the theory of Gr\"obner bases for bishuffle algebras, we consider the elements obtained from $e_k$ by adding a tail. For permutations $\sigma$ and $\tau$ let $\sigma\oplus\tau$ be their direct sum, the permutation whose matrix is block-diagonal with the blocks $\sigma$ and $\tau$. Put $e_k\oplus\tau=\sum_{\sigma\in\bS_k}\sgn(\sigma)\,\sigma\oplus\tau$, and for a set $\cW$ of permutations write $g\oplus\cW=\{g\oplus\tau:\tau\in\cW\}$.

\begin{introtheorem}[Theorem~\ref{thm:bwx}]\label{introthm:gb}
For every $k\ge1$ and every nonempty set $\cW$ of permutations, $e_k\oplus\cW$ is a Gr\"obner basis of the ideal it generates, both for the lexicographic ordering $<_\lex$ and for the opposite ordering $<_\lex^c$. The leading monomials are $k\cdots21\oplus\cW$ for $<_\lex$, and $12\cdots k\oplus\cW$ for $<_\lex^c$.
\end{introtheorem}

\begin{introcorollary}\label{introcor:bwx}
For every $k\ge1$, every nonempty set $\cW$ of permutations and every $N$,
\[
\#\Av_N\bigl(k\cdots21\oplus\cW\bigr)=\dim\bigl(\cS/(e_k\oplus\cW)\bigr)(N)=\#\Av_N\bigl(12\cdots k\oplus\cW\bigr).
\]
\end{introcorollary}

For a single tail $\cW=\{\tau\}$, the equality of the two numbers of avoiders is the theorem of Backelin, West and Xin \cite{BWX}, and for $k=3$ it is the theorem of Babson and West \cite{BW}. The extension to several tails is observed by Backelin, West and Xin right after their Proposition~2.3, and it is stated as a theorem by Burstein, Han, Kitaev and Zhang \cite[Theorem~2]{BHKZ}. So the enumerative statement is known. The point of Theorem~\ref{introthm:gb} is its algebraic form: the two sets of avoiders are two monomial bases of one and the same quotient of $\cS$.

The proof of Theorem~\ref{introthm:gb} follows the ideas of \cite{BW,BWX}. Starting from the set $\cW$, we split $\bS_N$ into classes. A point of a permutation is frozen if no copy of a tail from $\cW$ lies strictly south-east of it. The permutations with the same frozen points form a class $C$, and their remaining points run over all transversals of a Ferrers board $\lambda_C$. The quotient then splits along these classes,
\[
\bigl(\cS/(e_k\oplus\cW)\bigr)(N)\cong\bigoplus_{C}\bigl(R_{\lambda_C}/I_k(\lambda_C)\bigr)_{(1;1)},
\]
where $R_{\lambda_C}$ is the polynomial ring in the cells of $\lambda_C$, $I_k(\lambda_C)$ is the ladder determinantal ideal generated by its $k\times k$ minors, and $(1;1)$ denotes the multilinear part.

Thus everything reduces to ladder determinantal ideals, and this is where our approach differs from that of \cite{BWX}. The core technical result of Backelin, West and Xin is the shape-Wilf-equivalence of $12\cdots k$ and $k\cdots21$: for every Ferrers board $\lambda$, the number of transversals of $\lambda$ avoiding $12\cdots k$ in $\lambda$ equals the number of those avoiding $k\cdots21$ in $\lambda$. In \cite{BWX} this is proved combinatorially, by an explicit bijection. Here this step is replaced by a more geometric statement, which comes from the known Gr\"obner bases of matrix Schubert varieties. Namely, the ladder determinantal variety of $\lambda$, formed by the partial matrices on $\lambda$ whose north-west submatrices at the outer corners of $\lambda$ have rank less than $k$, is, up to an affine factor, the matrix Schubert variety of a vexillary permutation (Lemma~\ref{lem:board-vex}). By the theorems of Knutson and Miller \cite{KM} and of Knutson, Miller and Yong \cite{KMY}, this variety has two Gr\"obner degenerations, for an antidiagonal and for a diagonal term order, whose initial ideals are generated by the antidiagonal and by the diagonal terms of the $k\times k$ minors. The multilinear parts of the two degenerate rings have bases formed by the transversals of $\lambda$ without a decreasing $k$-chain and without an increasing $k$-chain in $\lambda$, respectively. Since a Gr\"obner degeneration preserves dimensions, the two numbers of transversals coincide (Theorem~\ref{thm:ladder}). So the shape-Wilf-equivalence of $12\cdots k$ and $k\cdots21$ appears as the comparison of two degenerations of one and the same determinantal variety.

\subsection{Organization of the paper}
Section~\ref{sec:algebras} recalls shuffle algebras and introduces bishuffle algebras, their Koszul duality, with the theorem of Carlitz, Scoville and Vaughan \cite{CSV} as an illustration, and commutative bishuffle algebras. Section~\ref{sec:examples} treats the main examples, namely the $\bS$-symmetric algebras and the algebras $\cB_V$, together with Schur--Weyl duality. Section~\ref{sec:noeth} is devoted to the Noetherianity questions. Section~\ref{sec:GB} develops Gr\"obner bases in $\cS$, with the quotients $\cS/(e_k)$ and $\cA_{1,1}$ as examples. Section~\ref{sec:ferrers} recalls Ferrers boards and deduces the shape-Wilf-equivalence of $12\cdots k$ and $k\cdots21$ from the results on matrix Schubert varieties. Section~\ref{sec:BWX} is devoted to the proof of Theorem~\ref{introthm:gb}.

\subsubsection*{{\bf {Acknowledgement}}}
The author is grateful to Doron Zeilberger for asking, after the author's talk at the Rutgers Experimental Mathematics Seminar in 2012, whether the theory of shuffle algebras and the Wilf-equivalence results established in~\cite{DK} could be generalized to classical pattern avoidance problems. That question directly motivated the present work.

\subsubsection*{{\bf AI usage statement}}
The project originated before the available AI tools were strong enough to contribute meaningfully. However, it would not have appeared without the help of Claude Opus 5.5, which was used extensively at the later stages for proofreading, discussions, numerical experiments and literature research.

\section{Shuffle and bishuffle algebras}\label{sec:algebras}

Throughout, $\kk$ is a field of characteristic zero and $[n]=\{1,\dots,n\}$. Permutations are written in one-line notation, and for a sequence $a_1\cdots a_k$ of distinct numbers $\st(a_1\cdots a_k)\in\bS_k$ denotes its standardization.

\subsection{Shuffle algebra and consecutive pattern avoidance}\label{sec:shuffle}

An \emph{$(m,n)$-shuffle} is a permutation $\sigma\in\bS_{m+n}$ with $\sigma(1)<\dots<\sigma(m)$ and $\sigma(m+1)<\dots<\sigma(m+n)$. The shuffles are the shortest representatives of the cosets in $\bS_{m+n}/(\bS_m\times\bS_n)$, and $\sigma$ corresponds to the decomposition $[m+n]=\sigma([m])\sqcup\sigma([m+1,m+n])$.

Let $\Ord$ be the category of finite totally ordered sets and order-preserving bijections. A \emph{collection} is a contravariant functor $\cM\colon\Ord\to\Vect$; it is determined by the spaces $\cM(n)=\cM([n])$. For $f\in\cM(n)$ and an ordered set $I$ with $|I|=n$ we write $f_I\in\cM(I)$ for the transport of $f$ along the unique order-preserving bijection $I\cong[n]$. The \emph{shuffle tensor product} of collections is
\[
(\mathscr P\boxtimes\mathscr Q)(I)=\bigoplus_{I=J\sqcup K}\mathscr P(J)\otimes\mathscr Q(K),
\]
the sum running over all decompositions of $I$ into two disjoint, not necessarily consecutive, subsets with the induced orders. An \emph{associative shuffle algebra} is a monoid for $\boxtimes$, that is, a collection $\cA$ with natural associative products $\cA(J)\otimes\cA(K)\to\cA(J\sqcup K)$; we always assume $\cA(\varnothing)=\kk$. On exponential generating series $f_\cM(t)=\sum_n\dim\cM(n)\,t^n/n!$ one has $f_{\mathscr P\boxtimes\mathscr Q}=f_{\mathscr P}f_{\mathscr Q}$.

Equivalently, an associative shuffle algebra is a sequence of spaces $\cA(n)$ with associative products
\[
\star_\sigma\colon\cA(m)\otimes\cA(n)\to\cA(m+n)
\]
indexed by the $(m,n)$-shuffles $\sigma$: the product $u\star_\sigma v$ is the image of $u_{\sigma([m])}\otimes v_{\sigma([m+1,m+n])}$ under the map $\cA(\sigma([m]))\otimes\cA(\sigma([m+1,m+n]))\to\cA([m+n])$.

The free shuffle algebra on one generator of arity one is the algebra $\AMR$ with $\AMR(I)=\kk\,\Sym(I)$, whose monomials are the linear arrangements of $I$ and whose product is concatenation of words \cite{Ronco,DK}. Thus $\AMR(n)=\kk\bS_n$ has the basis of permutations, and for $u\in\bS_m$ and $v\in\bS_n$ the product $u\star_\sigma v$ is the word $u$, written in the letters $\sigma([m])$, followed by the word $v$, written in the letters $\sigma([m+1,m+n])$. A monomial $p$ divides a monomial $w$ if some consecutive factor of $w$ standardizes to $p$; thus divisibility in $\AMR$ is containment of consecutive patterns. Consequently, for a set $P$ of permutations, the ideal $(P)$ is spanned by the permutations containing a consecutive pattern from $P$, and the monomial quotient $\AMR/(P)$ has a basis of the permutations avoiding all consecutive patterns from $P$. Monomial orderings, Gr\"obner bases, S-polynomials and the Buchberger criterion for shuffle algebras are developed in \cite{DK}, and the Anick-type resolution of $\AMR/(P)$ gives the cluster inversion formula of \cite[\S4]{DK}.

\subsection{Bishuffle algebras}\label{sec:bishuffle}

Let $\Ord^{(2)}$ be the category whose objects are pairs $(I,J)$ of finite ordered sets with $|I|=|J|$ and whose morphisms are pairs of order-preserving bijections. A \emph{bi-collection} is a contravariant functor $\cA\colon\Ord^{(2)}\to\Vect$, determined by the spaces $\cA(n)=\cA([n],[n])$.

\begin{definition}
\label{def::bishuffle}
The \emph{bishuffle tensor product} of bi-collections is
\[
(\mathscr P\boxtimes\mathscr Q)(I,J)=\bigoplus_{\substack{I=I_1\sqcup I_2,\ J=J_1\sqcup J_2\\ |I_1|=|J_1|}}\mathscr P(I_1,J_1)\otimes\mathscr Q(I_2,J_2).
\]
A \emph{bishuffle algebra} is a monoid for $\boxtimes$. We always assume $\cA(\varnothing,\varnothing)=\kk$.
\end{definition}

Equivalently, a bishuffle algebra is a sequence $\cA(n)$ with associative compositions 
$$\star_{\sigma,\tau}\colon\cA(n)\otimes\cA(m)\to\cA(n+m) \text{ indexed by pairs of $(n,m)$-shuffles.} $$
The natural generating series is $f_\cA(t)=\sum_n\dim\cA(n)\,t^n/(n!)^2$; it is multiplicative, since there are $\binom nk^2$ pairs of decompositions with $|I_1|=|J_1|=k$.

The free bishuffle algebra $F$ on one generator of arity one has as monomials on $(I,J)$ the words $(i_1,j_1)\cdots(i_n,j_n)$ with $I=\{i_1,\dots,i_n\}$ and $J=\{j_1,\dots,j_n\}$, that is, pairs of permutations $(\alpha,\beta)\in\Sym(I)\times\Sym(J)$. The product is concatenation of words, so $F(n)=\kk(\bS_n\times\bS_n)$ and $f_F(t)=1/(1-t)$. A monomial $p$ divides a monomial $w$ if a consecutive factor of the word of $w$ standardizes, in both coordinates separately, to $p$.

A \emph{monomial ordering} on $F$ is a family of well-orders on the sets of monomials on $(I,J)$, natural with respect to pairs of order-preserving bijections and such that every product map is strictly increasing in each argument. Examples are the lexicographic ordering of words of letters, with letters ordered by their first and then their second coordinate, and the \emph{$\beta$-first} ordering: $(\alpha,\beta)<(\alpha',\beta')$ if $\beta<\beta'$ lexicographically, or $\beta=\beta'$ and $\alpha<\alpha'$ lexicographically. Both are monomial orderings, because concatenation of words is an order embedding for the lexicographic order in each coordinate.

\subsection{Koszul theory}\label{sec:koszul}

Bishuffle algebras are the monoids in the monoidal category of bi-collections with the bishuffle tensor product. This category is $\kk$-linear and abelian, and $\boxtimes$ is exact in each argument. Hence the theory of quadratic algebras and of Koszul duality \cite{Priddy,PP,LV} applies to bishuffle algebras, exactly as it applies to shuffle algebras \cite{DK}. We only need algebras with one generator of arity one.

A bishuffle algebra is \emph{quadratic} if it is the quotient $A=F/(R)$ of the free bishuffle algebra $F$ of Section~\ref{sec:bishuffle} by the ideal generated by a subspace $R\subset F(2)$. Recall that $F(2)$ has the basis of the four monomials $(12,12)$, $(12,21)$, $(21,12)$ and $(21,21)$. The \emph{Koszul dual} of $A$ is
\[
A^!=F/(R^\perp),
\]
where $R^\perp\subset F(2)$ is the orthogonal complement of $R$ for the scalar product in which the four monomials are orthonormal; thus $(A^!)^!=A$. The algebra $A$ is \emph{Koszul} if $\operatorname{Tor}^A_i(\kk,\kk)$ is concentrated in arity $i$ for every $i$. The opposite of a monomial ordering of $F$ is again a monomial ordering. We use the following facts, whose proofs for associative algebras \cite{Priddy,PP,LV} carry over verbatim, as in the case of shuffle algebras \cite{DK}.
\begin{enumerate}
\item If $A$ is Koszul, then so is $A^!$, and $f_A(t)\,f_{A^!}(-t)=1$.
\item Suppose that $R$ is a Gr\"obner basis of $(R)$ for a monomial ordering of $F$, that is, the monomials of $F$ not divisible by any leading monomial of $R$ form a basis of $A$. Then $A$ is Koszul. Moreover, $R^\perp$ is a Gr\"obner basis of $(R^\perp)$ for the opposite ordering, and its leading monomials are the quadratic monomials that are not leading monomials of $R$.
\end{enumerate}
The identity in (1) is the Euler characteristic of the Koszul resolution of $\kk$, since $f$ is multiplicative for $\boxtimes$. Part (2) is the PBW criterion.

\subsubsection{Pairs of permutations without a common rise}
A \emph{common rise} of a pair $(\alpha,\beta)\in\bS_n\times\bS_n$ is a position $i\in[n-1]$ with $\alpha(i)<\alpha(i+1)$ and $\beta(i)<\beta(i+1)$. Let $\omega_n$ be the number of pairs without a common rise; for $n=0,1,\dots,6$ these numbers are $1,1,3,19,211,3651,90921$.

\begin{proposition}\label{prop:csv}
Let $A=F/(R)$, where $R\subset F(2)$ is spanned by $(12,12)$, $(21,21)$ and $(12,21)-(21,12)$.
\begin{enumerate}
\item For all $n$, $\dim A(n)=1$, and $A$ is Koszul.
\item The Koszul dual $A^!=F/\bigl((12,21)+(21,12)\bigr)$ has a basis of the pairs $(\alpha,\beta)\in\bS_n\times\bS_n$ with no position $i$ at which $\alpha$ falls and $\beta$ rises. In particular $\dim A^!(n)=\omega_n$.
\item Consequently
\[
\sum_{n\ge0}\omega_n\frac{t^n}{(n!)^2}=\Bigl(\sum_{n\ge0}(-1)^n\frac{t^n}{(n!)^2}\Bigr)^{-1},
\]
which is the theorem of Carlitz, Scoville and Vaughan \cite{CSV}.
\end{enumerate}
\end{proposition}

\begin{proof}
First, we show that $R$ is a Gr\"obner basis for the $\beta$-first ordering of Section~\ref{sec:bishuffle}. The leading monomials of the three relations are $(12,12)$, $(21,21)$ and $(12,21)$, since $(12,21)$ has the larger second coordinate. So a monomial is normal if and only if every two adjacent letters form the pattern $(21,12)$, that is, if and only if its first coordinate is decreasing and its second coordinate is increasing. The only normal monomial of arity $n$ is $(n\cdots21,\,12\cdots n)$, and it remains to show that it is nonzero in $A$. Let $\varphi\colon F(n)\to\kk$ send a word to $1$ if all its letters $(a,b)$ satisfy $a+b=n+1$, and to $0$ otherwise. Transported to a pair of two-element sets $(\{i_1<i_2\},\{j_1<j_2\})$, the monomial relations become the factors $(i_1,j_1)(i_2,j_2)$ and $(i_2,j_2)(i_1,j_1)$, whose two letters have different sums $a+b$; so $\varphi$ vanishes on the words containing them. The binomial relation becomes $(i_1,j_2)(i_2,j_1)-(i_2,j_1)(i_1,j_2)$, the exchange of two adjacent letters, and $\varphi$ only depends on the set of letters. Hence $\varphi$ vanishes on the ideal $(R)$, while $\varphi(n\cdots21,\,12\cdots n)=1$. So $\dim A(n)=1$, $R$ is a Gr\"obner basis, and $A$ is Koszul by fact~(2).

Second, $R^\perp$ is spanned by $(12,21)+(21,12)$. For the opposite of the $\beta$-first ordering its leading monomial is $(21,12)$, the quadratic monomial that is not a leading monomial of $R$, and by fact~(2) it forms a Gr\"obner basis. So $A^!(n)$ has a basis of the monomials without a consecutive factor of the pattern $(21,12)$, that is, of the pairs $(\alpha,\beta)$ with no position $i$ such that $\alpha(i)>\alpha(i+1)$ and $\beta(i)<\beta(i+1)$. Replacing $\alpha$ by its complement $i\mapsto n+1-\alpha(i)$ turns the falls of $\alpha$ into rises, so these pairs are in bijection with the pairs without a common rise.

Finally, by fact~(1), $f_{A^!}(t)=1/f_A(-t)$, and $f_A(t)=\sum_nt^n/(n!)^2$ by the first part.
\end{proof}

\subsection{Commutative bishuffle algebras}\label{sec:commutative}

The \emph{opposite} of an $(m,n)$-shuffle \\ $\sigma=(\sigma(1)<\dots<\sigma(m);\ \sigma(m+1)<\dots<\sigma(m+n))$ 
\[
\text{is the $(n,m)$-shuffle }
\sigma^{op}=\bigl(\sigma(m+1)<\dots<\sigma(m+n);\ \sigma(1)<\dots<\sigma(m)\bigr),
\]
that is, $\sigma^{op}(i)=\sigma(m+i)$ for $i\le n$ and $\sigma^{op}(n+j)=\sigma(j)$ for $j\le m$. It corresponds to the same decomposition of $[m+n]$, with the two parts exchanged.

\begin{definition}\label{def:commutative}
A bishuffle algebra $\cA$ is \emph{commutative} if
\[
u\star_{\sigma,\tau}v=v\star_{\sigma^{op},\tau^{op}}u
\]
for all $u\in\cA(m)$, $v\in\cA(n)$ and all pairs of $(m,n)$-shuffles $\sigma,\tau$.
\end{definition}

The products $\star_{\sigma,\tau}$ are the components of the product $\cA\boxtimes\cA\to\cA$, and the pair $(\sigma^{op},\tau^{op})$ indexes the same decompositions of the values and of the positions as $(\sigma,\tau)$, with the two factors exchanged. So $\cA$ is commutative if and only if its product is invariant under the exchange of the two factors of $\boxtimes$, that is, if $\cA$ is a commutative monoid for $\boxtimes$.

Consider the relations in arity two
\[
r_1=(21,21)-(12,12),\qquad r_2=(12,21)-(21,12),
\]
and put $\cS=F/(r_1,r_2)$. Transported to $(\{i_1<i_2\},\{j_1<j_2\})$, the relations read $(i_2,j_2)(i_1,j_1)=(i_1,j_1)(i_2,j_2)$ and $(i_1,j_2)(i_2,j_1)=(i_2,j_1)(i_1,j_2)$, so together they allow us to exchange any two adjacent letters.

\begin{proposition}\label{prop:S}
\begin{enumerate}
\item Sending a word to its set of letters identifies $\cS(I,J)$ with $\kk\Bij(J,I)$; under this identification the product of two bijections is their disjoint union.
\item In particular $\cS(n)\cong\kk\bS_n$, where the class of $(\alpha,\beta)$ corresponds to the permutation $\alpha\beta^{-1}$, and $f_\cS(t)=\sum_{n=0}^{\infty} \frac{n!}{n!^2}t^n= e^t$.
\item $\{r_1,r_2\}$ is a Gr\"obner basis of $(r_1,r_2)$ for every monomial ordering of $F$.
\end{enumerate}
\end{proposition}

\begin{proof}
Two words are equal in $\cS$ if and only if they are connected by exchanges of adjacent letters, that is, if and only if they have the same set of letters; and the set of letters of $(\alpha,\beta)$ is the graph of the bijection $\beta(k)\mapsto\alpha(k)$. This gives (1) and (2). For (3), a monomial ordering chooses one of the two monomials in each $r_i$ as leading. In each of the four possible cases, the two leading monomials share the same pattern ($12$ or $21$) in one coordinate and exhaust both patterns in the other. Hence the normal monomials are the pairs whose distinguished coordinate is $\id$, respectively $w_0$. There are $n!=\dim\cS(n)$ of them, and normal monomials always span the quotient.
\end{proof}

\subsubsection{Bishuffle products of permutations}
From now on we identify $\cS(n)$ with $\kk\bS_n$ and work directly with permutations.

Let us describe the products of $\cS$ on permutations in several different ways. Permutations are maps from positions to values, multiplied as maps: $(uv)(i)=u(v(i))$. For $u\in\bS_n$ and $v\in\bS_m$ we write $u\oplus v=u(1)\cdots u(n)\,(v(1)+n)\cdots(v(m)+n)\in\bS_{n+m}$ for their direct sum. For $L\subset[N]$ with $|L|=k$ we write $\pi_L$ for the $(k,N-k)$-shuffle with $\pi_L([k])=L$. By Proposition~\ref{prop:S}(1), in an index $(I,J)$ of the bi-collection $\cS$ the set $I$ consists of values and the set $J$ of positions. So a product $\star_{\sigma,\tau}$ of $\cS$ is indexed by a shuffle $\sigma$ of the values and a shuffle $\tau$ of the positions.

\begin{proposition}\label{prop:products}
Let $u\in\bS_n$ and $v\in\bS_m$, and let $\sigma$ and $\tau$ be $(n,m)$-shuffles. Then
\[
u\star_{\sigma,\tau}v=\sigma\,(u\oplus v)\,\tau^{-1}.
\]
Thus $w=u\star_{\sigma,\tau}v$ is the permutation of $[n+m]$ whose positions $\tau([n])$ carry the values $\sigma([n])$ in the pattern $u$, and whose positions $\tau([n+1,n+m])$ carry the values $\sigma([n+1,n+m])$ in the pattern $v$. Conversely, for $w\in\bS_{n+m}$ and a set $L$ of $n$ positions, $w=u\star_{\sigma,\pi_L}v$ holds if and only if $\sigma=\pi_{w(L)}$, and $u$ and $v$ are the patterns formed by $w$ on $L$ and on $[n+m]\setminus L$.
\end{proposition}

\begin{proof}
By Proposition~\ref{prop:S}(1), the product of $\cS$ is the disjoint union of bijections from positions to values. Since $\sigma$ and $\tau$ are increasing on $[n]$, the transport of $u$ to the values $\sigma([n])$ and the positions $\tau([n])$ is the restriction of $\sigma u\tau^{-1}$ to $\tau([n])$; similarly, the transport of $v$ is the restriction of $\sigma(u\oplus v)\tau^{-1}$ to $\tau([n+1,n+m])$. Their union is $\sigma(u\oplus v)\tau^{-1}$. Equivalently, in $F$ the product of the words $(\alpha,\beta)$ and $(\alpha',\beta')$ along $(\sigma,\tau)$ is the word $\bigl(\sigma(\alpha\oplus\alpha'),\tau(\beta\oplus\beta')\bigr)$, and the map $(\alpha,\beta)\mapsto\alpha\beta^{-1}$ of Proposition~\ref{prop:S}(2) sends it to $\sigma(\alpha\beta^{-1}\oplus\alpha'\beta'^{-1})\tau^{-1}$. For the converse, the positions $L=\tau([n])$ determine $\tau=\pi_L$, the values $w(L)=\sigma([n])$ determine $\sigma$, and then $u$ and $v$ are the patterns of $w$ on $L$ and on its complement.
\end{proof}

\begin{remark}\label{rem:products}
The factor $\sigma(u\oplus v)$ is the product $u\star_\sigma v$ of the shuffle algebra $\AMR$ of Section~\ref{sec:shuffle}: it is the concatenation of the word $u$, written in the values $\sigma([n])$, with the word $v$, written in the values $\sigma([n+1,n+m])$. Composing with $\tau^{-1}$ on the right moves the first $n$ letters of this word to the positions $\tau([n])$ and the last $m$ letters to the positions $\tau([n+1,n+m])$, keeping their order. Hence
\[
u\star_{\sigma,\tau}v=(u\star_\sigma v)\,\tau^{-1},
\]
a product of the shuffle algebra in the values followed by a shuffle of the positions. Symmetrically, $(u\oplus v)\tau^{-1}$ interleaves the word $u$ with the word $v$ shifted by $n$, and $u\star_{\sigma,\tau}v$ is obtained from it by relabelling the values by $\sigma$.
\end{remark}

\begin{example}\label{ex:products}
\begin{enumerate}
\item For $n=m=1$ both shuffles are $12$ or $21$, and $1\star_{\sigma,\tau}1=\sigma\tau^{-1}$. In the free algebra $F$ the four products of the generator with itself are the four words $(\sigma,\tau)$. In $\cS$ the products for $(\sigma,\tau)=(12,12)$ and $(21,21)$ both equal $12$, and those for $(12,21)$ and $(21,12)$ both equal $21$. These two coincidences are the relations $r_1$ and $r_2$.
\item Let $u=21$ and $v=1$. The $(2,1)$-shuffles are $123$, $132$ and $231$, and $u\oplus v=213$. The nine products $21\star_{\sigma,\tau}1$ are
\[
\begin{array}{c|ccc}
\sigma\ \backslash\ \tau & 123 & 132 & 231\\\hline
123 & 213 & 231 & 321\\
132 & 312 & 321 & 231\\
231 & 321 & 312 & 132
\end{array}
\]
For instance, for $\sigma=132$ and $\tau=231$ the shuffle product is $21\star_{132}1=312$, the word $21$ written in the values $\{1,3\}$ followed by the value $2$; moving its first two letters to the positions $\{2,3\}$ gives $231$. Every $w\in\bS_3$ appears as many times as it has inversions: by the last part of Proposition~\ref{prop:products}, a pair $(\sigma,\tau)$ with $21\star_{\sigma,\tau}1=w$ is the same as an occurrence of $21$ in $w$.
\item For $u=132$, $v=21$, $\sigma=\pi_{\{1,3,5\}}=13524$ and $\tau=\pi_{\{2,3,5\}}=23514$ one gets $u\star_{\sigma,\tau}v=41523$: the positions $2,3,5$ carry the values $1,5,3$ in the pattern $132$, and the positions $1,4$ carry the values $4,2$ in the pattern $21$.
\end{enumerate}
\end{example}

\begin{corollary}\label{cor:S-commutative}
The bishuffle algebra $\cS$ is commutative in the sense of Definition~\ref{def:commutative}. Moreover, it is the free commutative bishuffle algebra on one generator of arity one: every commutative quotient of $F$ is a quotient of $\cS$.
\end{corollary}

\begin{proof}
First, let $u\in\bS_m$ and $v\in\bS_n$, and let $c\in\bS_{m+n}$ be given by $c(i)=m+i$ for $i\le n$ and $c(n+j)=j$ for $j\le m$. Then $\sigma^{op}=\sigma c$ and $\tau^{op}=\tau c$ for all $(m,n)$-shuffles $\sigma$ and $\tau$. Moreover $c\,(v\oplus u)\,c^{-1}=u\oplus v$. Indeed, for $j\le m$ we have $c^{-1}(j)=n+j$ and $c\bigl((v\oplus u)(n+j)\bigr)=c(n+u(j))=u(j)$, and for $i\le n$ we have $c^{-1}(m+i)=i$ and $c\bigl((v\oplus u)(i)\bigr)=c(v(i))=m+v(i)$. Hence, by Proposition~\ref{prop:products},
\[
v\star_{\sigma^{op},\tau^{op}}u=\sigma c\,(v\oplus u)\,c^{-1}\tau^{-1}=\sigma\,(u\oplus v)\,\tau^{-1}=u\star_{\sigma,\tau}v .
\]

Second, by Example~\ref{ex:products}(1), the four products of the generator $x$ of $F$ with itself are the words $(\sigma,\tau)$ with $\sigma,\tau\in\{12,21\}$, and the opposite of the $(1,1)$-shuffle $12$ is $21$. So the relations $x\star_{\sigma,\tau}x=x\star_{\sigma^{op},\tau^{op}}x$ are, up to sign, $r_1$ and $r_2$. They hold in every commutative quotient of $F$, so every such quotient is a quotient of $\cS$, which is commutative by the first part.
\end{proof}

\subsubsection{Koszulness and the Koszul dual}
By construction, $\cS=F/(r_1,r_2)$ is a quadratic bishuffle algebra, and $\{r_1,r_2\}$ is a Gr\"obner basis by Proposition~\ref{prop:S}(3). Hence $\cS$ is Koszul by fact~(2) of Section~\ref{sec:koszul}. Its Koszul dual is defined by the anticommutativity relations,
\[
\cS^!=F\big/\bigl((12,12)+(21,21),\ (12,21)+(21,12)\bigr).
\]
We call a bishuffle algebra $\cA$ \emph{skew-commutative} if
\[
u\star_{\sigma,\tau}v=(-1)^{mn}\,v\star_{\sigma^{op},\tau^{op}}u
\]
for all $u\in\cA(m)$, $v\in\cA(n)$ and all pairs of $(m,n)$-shuffles $\sigma,\tau$.

\begin{proposition}\label{prop:S-dual}
\begin{enumerate}
\item Sending a word $(\alpha,\beta)$ to $\sgn(\beta)\,\alpha\beta^{-1}$ identifies $\cS^!(n)$ with $\kk\bS_n$. In particular, $\cS^!$ has a basis of permutations, $\dim\cS^!(n)=n!$ and $f_{\cS^!}(t)=e^t$.
\item In this basis the product of $u\in\bS_m$ and $v\in\bS_n$ is
\[
u\star_{\sigma,\tau}v=\sgn(\tau)\,\sigma\,(u\oplus v)\,\tau^{-1}.
\]
\item $\cS^!$ is the free skew-commutative bishuffle algebra on one generator of arity one: it is skew-commutative, and every skew-commutative quotient of $F$ is a quotient of $\cS^!$.
\end{enumerate}
\end{proposition}

\begin{proof}
First, let $\psi\colon F(n)\to\kk\bS_n$ be the linear map with $\psi(\alpha,\beta)=\sgn(\beta)\,\alpha\beta^{-1}$. Transported to $(\{i_1<i_2\},\{j_1<j_2\})$, the defining relations of $\cS^!$ read $(i_1,j_1)(i_2,j_2)+(i_2,j_2)(i_1,j_1)$ and $(i_1,j_2)(i_2,j_1)+(i_2,j_1)(i_1,j_2)$. So the ideal they generate is spanned by the sums $x+x'$, where the word $x'$ is obtained from the word $x$ by exchanging two adjacent letters. Such an exchange does not change $\alpha\beta^{-1}$ and multiplies $\sgn(\beta)$ by $-1$, so $\psi$ vanishes on this ideal and induces a surjection $\cS^!(n)\to\kk\bS_n$. On the other hand, modulo the relations every word is, up to sign, equal to the word with the same letters arranged by increasing second coordinates, so $\dim\cS^!(n)\le n!$. This proves (1). The series is $\sum_nn!\,t^n/(n!)^2=e^t$, in accordance with $f_\cS(t)f_{\cS^!}(-t)=1$.

Second, the product of the words $(\alpha,\beta)$ and $(\alpha',\beta')$ along $(\sigma,\tau)$ is the word $\bigl(\sigma(\alpha\oplus\alpha'),\tau(\beta\oplus\beta')\bigr)$, and $\psi$ sends it to
\[
\sgn(\tau)\sgn(\beta)\sgn(\beta')\,\sigma\bigl(\alpha\beta^{-1}\oplus\alpha'\beta'^{-1}\bigr)\tau^{-1}.
\]
This proves (2).

Finally, let $c\in\bS_{m+n}$ be the permutation from the proof of Corollary~\ref{cor:S-commutative}, so that $\sigma^{op}=\sigma c$, $\tau^{op}=\tau c$ and $c\,(v\oplus u)\,c^{-1}=u\oplus v$. It has $mn$ inversions, so by (2)
\[
v\star_{\sigma^{op},\tau^{op}}u=(-1)^{mn}\sgn(\tau)\,\sigma c\,(v\oplus u)\,c^{-1}\tau^{-1}=(-1)^{mn}\,u\star_{\sigma,\tau}v .
\]
For $m=n=1$ the skew-commutativity relations $x\star_{\sigma,\tau}x+x\star_{\sigma^{op},\tau^{op}}x$ of the generator $x$ of $F$ are the defining relations of $\cS^!$, so every skew-commutative quotient of $F$ is a quotient of $\cS^!$.
\end{proof}

\subsubsection{Permutation matrices}
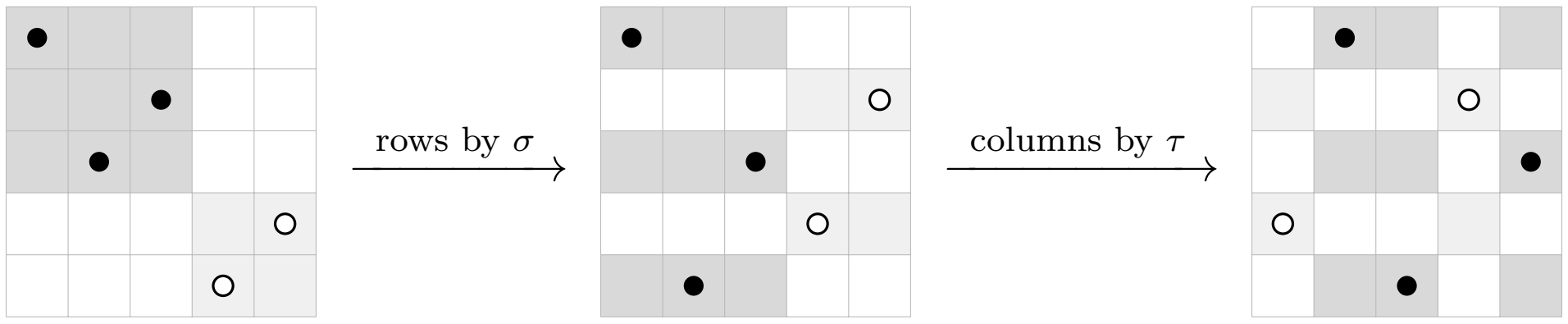
\begin{figure}[ht]
\centering
\begin{tikzpicture}[scale=0.5]
\foreach \r in {1,2,3}{\foreach \c in {1,2,3}{\shadecell[gray!30]{\r}{\c}}}
\foreach \r in {4,5}{\foreach \c in {4,5}{\shadecell[gray!12]{\r}{\c}}}
\gridlines{5}{5}
\pt{1}{1}\pt{3}{2}\pt{2}{3}\opt{5}{4}\opt{4}{5}
\end{tikzpicture}
\raisebox{1.1cm}{$\ \xrightarrow{\ \text{rows by }\sigma\ }\ $}
\begin{tikzpicture}[scale=0.5]
\foreach \r in {1,3,5}{\foreach \c in {1,2,3}{\shadecell[gray!30]{\r}{\c}}}
\foreach \r in {2,4}{\foreach \c in {4,5}{\shadecell[gray!12]{\r}{\c}}}
\gridlines{5}{5}
\pt{1}{1}\pt{5}{2}\pt{3}{3}\opt{4}{4}\opt{2}{5}
\end{tikzpicture}
\raisebox{1.1cm}{$\ \xrightarrow{\ \text{columns by }\tau\ }\ $}
\begin{tikzpicture}[scale=0.5]
\foreach \r in {1,3,5}{\foreach \c in {2,3,5}{\shadecell[gray!30]{\r}{\c}}}
\foreach \r in {2,4}{\foreach \c in {1,4}{\shadecell[gray!12]{\r}{\c}}}
\gridlines{5}{5}
\pt{1}{2}\pt{5}{3}\pt{3}{5}\opt{4}{1}\opt{2}{4}
\end{tikzpicture}
\caption{The product $132\star_{\sigma,\tau}21=41523$ of Example~\ref{ex:products}(3) in terms of permutation matrices, for $\sigma=13524$ and $\tau=23514$. Left: the block-diagonal matrix $P_{132\oplus21}=P_{13254}$; the entries of the block $P_{132}$ are black and those of the block $P_{21}$ white. Middle: the row $i$ moved to the row $\sigma(i)$. Right: the column $j$ moved to the column $\tau(j)$, which gives $P_{41523}$. The shaded cells are the images of the two diagonal blocks.}\label{fig:matrices}
\end{figure}

For $u\in\bS_m$ let $P_u\in\Mat_m$ be its permutation matrix, whose entry in the row $u(j)$ and the column $j$ is $1$ for every $j\in[m]$ and whose other entries are $0$. Thus $P_ue_j=e_{u(j)}$ and $P_{uv}=P_uP_v$; the rows of $P_u$ are indexed by the values and its columns by the positions. Let $v\in\bS_n$, and let $\sigma,\tau$ be $(m,n)$-shuffles. By Proposition~\ref{prop:products},
\[
P_{u\star_{\sigma,\tau}v}=P_\sigma\begin{pmatrix}P_u&0\\0&P_v\end{pmatrix}P_\tau^{-1}.
\]
In other words, $u\star_{\id,\id}v=u\oplus v$ is the block-diagonal matrix with the blocks $P_u$ and $P_v$, and $u\star_{\sigma,\tau}v$ is obtained from it by moving its $i$-th row to the row $\sigma(i)$ and its $j$-th column to the column $\tau(j)$, for all $i,j\in[m+n]$. Since $\sigma$ and $\tau$ are shuffles, the rows of each block keep their relative order, and so do its columns: the block $P_u$ is spread over the rows $\sigma([m])$ and the columns $\tau([m])$, and the block $P_v$ over the rows $\sigma([m+1,m+n])$ and the columns $\tau([m+1,m+n])$. Figure~\ref{fig:matrices} shows the product of Example~\ref{ex:products}(3).

So the free commutative bishuffle algebra $\cS$ of Corollary~\ref{cor:S-commutative} is spanned by the permutation matrices, its generator being the $1\times1$ matrix $(1)$, and its products interleave the rows and, independently, the columns of the two diagonal blocks. Commutativity says that the matrix does not change if the two blocks are listed in the opposite order and both shuffles are replaced by their opposites. We use these pictures throughout the paper: a permutation $u$ is drawn as the matrix of the linear map $e_j\mapsto e_{u(j)}$, with the columns indexed by the positions and the rows by the values.

\subsubsection{Patterns and substitutions}
For $w\in\bS_n$ and a $k$-subset $L=\{l_1<\dots<l_k\}\subset[n]$ of positions we write $w|_L=\st(w(l_1)\cdots w(l_k))\in\bS_k$. We say that $L$ is an \emph{occurrence} of $p\in\bS_k$ in $w$ if $w|_L=p$; thus $w$ contains $p$ as a classical pattern if and only if $p$ has an occurrence in $w$. For $\sigma\in\bS_k$ let $w_{L\leftarrow\sigma}\in\bS_n$ be the permutation that coincides with $w$ outside $L$ and takes on $L$ the values $w(L)$ arranged in the pattern $\sigma$. For $g=\sum_\sigma c_\sigma\sigma\in\kk\bS_k$ put
\[
g[w,L]=\sum_{\sigma\in\bS_k}c_\sigma\,w_{L\leftarrow\sigma}\in\kk\bS_n .
\]
By Proposition~\ref{prop:products}, $w_{L\leftarrow\sigma}=\sigma\star_{\pi_{w(L)},\pi_L}w|_{[n]\setminus L}$, and hence
\[
g[w,L]=g\star_{\pi_{w(L)},\pi_L}w|_{[n]\setminus L}:
\]
substituting $g$ into the positions $L$ of $w$ is multiplying $g$ by the pattern of $w$ on the complement of $L$. Hence the ideal of $\cS$ generated by a set $G$ of elements is
\[
(G)(n)=\operatorname{span}\bigl\{\,g[w,L]\ :\ g\in G,\ w\in\bS_n,\ L\subset[n],\ |L|=\text{arity of }g\,\bigr\}.
\]

\begin{corollary}\label{cor:classical}
For a set $P$ of permutations, the monomial quotient $\cS/(P)$ has a basis of the permutations avoiding all patterns from $P$, and
\[
f_{\cS/(P)}(t)=\sum_{n\ge0}\frac{\#\Av_n(P)}{(n!)^2}\,t^n .
\]
\end{corollary}

\begin{proof}
By the formula above, $(P)(n)$ is spanned by the permutations $p[w,L]=w$, where $L$ is an occurrence of some $p\in P$ in $w$.
\end{proof}

\begin{remark}\label{rem:patterns}
Thus classical pattern avoidance is the theory of monomial ideals of the one-generated free commutative bishuffle algebra $\cS$, exactly as consecutive pattern avoidance is the theory of monomial ideals of $\AMR$ \cite{DK}. The classes $\Av(P)$ are the monomial bases of the quotients $\cS/(P)$, the numbers $\#\Av_n(P)$ are the coefficients of their Hilbert series $f_{\cS/(P)}(t)$, and two sets of patterns are Wilf-equivalent if and only if the corresponding monomial quotients have the same Hilbert series. The simplest example is the algebra $A$ of Proposition~\ref{prop:csv}, which is the quotient $\cS/(12)$. Indeed, $r_1$ and $r_2$ lie in $(R)$, and modulo $r_1$ the monomials $(12,12)$ and $(21,21)$ both represent the permutation $12$. So the first part of the proof of Proposition~\ref{prop:csv} is the equality $\dim(\cS/(12))(n)=\#\Av_n(12)=1$ of Corollary~\ref{cor:classical}, and Proposition~\ref{prop:csv} identifies the Koszul dual of the simplest pattern-avoidance algebra with an algebra counting pairs of permutations.
\end{remark}

\begin{remark}\label{rem:twisted}
A collection $\cM$ is \emph{symmetric} if the functor $\cM$ extends from ordered sets to the category of finite sets and all bijections, that is, if every $\cM(n)$ is a representation of $\bS_n$. Symmetric collections are the linear species of Joyal \cite{Joyal,AM}, and on them the shuffle tensor product of Section~\ref{sec:shuffle} is the Cauchy product. A shuffle algebra on a symmetric collection whose products are natural for all bijections is known as a \emph{twisted associative algebra} \cite{Dotsenko}, and a commutative one as a \emph{twisted commutative algebra}. Twisted Lie algebras go back to Barratt \cite{Barratt}, and twisted commutative algebras were studied systematically by Sam and Snowden \cite{SS-tca,SS-GL,SS}. In characteristic zero, Schur--Weyl duality identifies twisted commutative algebras with commutative algebras carrying a polynomial action of $GL_\infty$ \cite{SS-tca}. For instance, the modules over the twisted commutative algebra with $\cM(n)=\kk$ for all $n$ are the FI-modules of Church, Ellenberg and Farb \cite{CEF,SS-tca}, and Noetherianity is one of the central questions of the subject \cite{SS-tca,NSS}.

We extend this language to bishuffle algebras in order to have an algebraic framework for pattern avoidance. The $\bS$-symmetric commutative bishuffle algebras of Section~\ref{sec:symdef} are twisted commutative algebras in two sets of variables, that is, commutative algebras with a polynomial action of $GL_\infty\times GL_\infty$ \cite{SS-tca}. The free algebra $\cS$ corresponds to the polynomial ring in the entries of an infinite matrix, and the description of symmetric ideals by upper sets of Young's lattice in Proposition~\ref{prop:symideals} corresponds to the classification of the $GL\times GL$-stable ideals of polynomial rings on matrices by De Concini, Eisenbud and Procesi \cite{DCEP}; for a recent classification of this kind for other twisted algebras see \cite{DL}. However, the ideals describing pattern avoidance are the monomial ideals of $\cS$, which in general are not stable under permutations. They are not visible in the language of twisted algebras, and this is why we keep only the order-preserving bijections.
\end{remark}

\section{Symmetric commutative bishuffle algebras}\label{sec:examples}

In this section we describe the main examples of commutative bishuffle algebras. Section~\ref{sec:symdef} treats the quotients of $\cS$ by ideals that are stable under all permutations of positions and of values. In every arity such a quotient is the quotient of a group algebra by a two-sided ideal, so its dimension is given by representation theory. Section~\ref{sec:End} treats the algebra $\cB_V$ with the components $\End(V^{\otimes n})$ and identifies its invariants with $\bS$-symmetric algebras by Schur--Weyl duality. Section~\ref{sec:noeth} discusses the Noetherianity of these algebras.

\subsection{\texorpdfstring{$\bS$}{S}-symmetric commutative bishuffle algebras}\label{sec:symdef}

Recall from Section~\ref{sec:commutative} that permutations are maps from positions to values, multiplied as maps. So left multiplication by $\bS_N$ permutes the values, and right multiplication permutes the positions. Thus every component $\cS(N)=\kk\bS_N$ of the free algebra is an $\bS_N\times\bS_N$-module, and by Proposition~\ref{prop:products} the product $u\star_{\sigma,\tau}v=\sigma\,(u\oplus v)\,\tau^{-1}$ is the image of $u\oplus v$ under the left multiplication by $\sigma$ and the right multiplication by $\tau^{-1}$.

\begin{definition}\label{def:symmetric}
An ideal $\cJ\subset\cS$ is \emph{symmetric} if every component $\cJ(N)$ is a two-sided ideal of the group algebra $\kk\bS_N$, that is, if it is stable under the permutations of positions and of values. The quotient $\cS/\cJ$ by a symmetric ideal is called an \emph{$\bS$-symmetric commutative bishuffle algebra}. Its component of arity $N$ is the quotient of $\kk\bS_N$ by the two-sided ideal $\cJ(N)$.
\end{definition}

For a partition $\lambda$ of $N$ let $V_\lambda$ be the irreducible representation of $\bS_N$, let $f^\lambda=\dim V_\lambda$ be the number of standard Young tableaux of shape $\lambda$, and let $\End(V_\lambda)$ be the corresponding block in the Wedderburn decomposition $\kk\bS_N=\bigoplus_{\lambda\vdash N}\End(V_\lambda)$. As an $\bS_N\times\bS_N$-module, each block $\End(V_\lambda)\cong V_\lambda\otimes V_\lambda^*$ is irreducible. We write $\lambda'$ for the conjugate partition and $\ell(\lambda)$ for the number of rows of $\lambda$; recall that $f^{\lambda'}=f^\lambda$.

\begin{lemma}\label{lem:young}
Let $g\in\kk\bS_k$, $w\in\bS_N$ and $L\subset[N]$ with $|L|=k$. Then $g[w,L]=u\,g\,v$ for some $u,v\in\bS_N$, where $\kk\bS_k\subset\kk\bS_N$ via the first $k$ letters. If $J\subset\kk\bS_k$ is a two-sided ideal, then the elements $h[w,L]$, for $h\in J$ and all $w$ and $L$, span the two-sided ideal of $\kk\bS_N$ generated by $J$.
\end{lemma}

\begin{proof}
Let $M=w(L)$, and let $\rho\in\bS_{N-k}$ be the pattern of $w$ on $[N]\setminus L$. By Proposition~\ref{prop:products}, $w_{L\leftarrow\sigma}=\sigma\star_{\pi_M,\pi_L}\rho=\pi_M\,(\sigma\oplus\rho)\,\pi_L^{-1}$ for every $\sigma\in\bS_k$, so $g[w,L]=\pi_M\,g\,(\id_k\oplus\rho)\,\pi_L^{-1}$. Conversely, every $u\in\bS_N$ can be written as $u=\pi_M(\alpha\oplus\beta)$ with $M=u([k])$, and every $v\in\bS_N$ as $v=(\gamma\oplus\delta)\pi_L^{-1}$ with $L=v^{-1}([k])$, where $\alpha,\gamma\in\bS_k$ and $\beta,\delta\in\bS_{N-k}$. For $h\in J$ this gives
\[
u\,h\,v=\pi_M\,(\alpha h\gamma)\,(\id_k\oplus\beta\delta)\,\pi_L^{-1}=(\alpha h\gamma)[w',L],
\]
where $w'=\pi_M(\id_k\oplus\beta\delta)\pi_L^{-1}$, and $\alpha h\gamma\in J$.
\end{proof}

\begin{proposition}\label{prop:symideals}
For a symmetric ideal $\cJ\subset\cS$ let $\Lambda(\cJ)$ be the set of partitions $\lambda$ such that $\End(V_\lambda)\subset\cJ(|\lambda|)$. Then $\Lambda(\cJ)$ is an upper set of Young's lattice, and $\cJ\mapsto\Lambda(\cJ)$ is a bijection from symmetric ideals onto upper sets, whose inverse sends $\Lambda$ to the ideal $\cJ_\Lambda$ with
\[
\cJ_\Lambda(N)=\bigoplus_{\lambda\vdash N,\ \lambda\in\Lambda}\End(V_\lambda),\qquad \dim(\cS/\cJ_\Lambda)(N)=\sum_{\lambda\vdash N,\ \lambda\notin\Lambda}(f^\lambda)^2 .
\]
The symmetric ideal generated by the block $\End(V_\mu)$ corresponds to $\{\lambda:\lambda\supseteq\mu\}$.
\end{proposition}

\begin{proof}
Since $\kk$ has characteristic zero, $\kk\bS_N$ is semisimple, and its two-sided ideals are the sums of blocks. By Lemma~\ref{lem:young}, the ideal of $\cS$ generated by $\End(V_\mu)\subset\kk\bS_k$ is, in arity $N$, the two-sided ideal of $\kk\bS_N$ generated by $\End(V_\mu)$. This is the sum of the blocks $\End(V_\lambda)$ on which $\End(V_\mu)$ acts nontrivially, that is, for which the restriction of $V_\lambda$ to $\bS_k$ contains $V_\mu$. By the branching rule this happens if and only if $\mu\subseteq\lambda$. Hence $\Lambda(\cJ)$ is an upper set, and $\cJ_\Lambda$ is an ideal for every upper set $\Lambda$.
\end{proof}

For example, the antisymmetrizer $e_n=\sum_{\sigma\in\bS_n}\sgn(\sigma)\,\sigma$ and the symmetrizer $s_n=\sum_{\sigma\in\bS_n}\sigma$ are $e_n=n!\,\epsilon_{(1^n)}$ and $s_n=n!\,\epsilon_{(n)}$, where $\epsilon_\lambda$ is the central idempotent of the block of $\lambda$. So $(e_n)$, $(s_m)$ and $(e_n,s_m)$ are symmetric ideals, corresponding to the upper sets $\{\lambda:\ell(\lambda)\ge n\}$, $\{\lambda:\lambda_1\ge m\}$ and their union. On the other hand, not every ideal is symmetric; for instance, the component of arity three of the ideal generated by $123-213$ is spanned by this element alone.

\subsubsection{The sign twist and the symmetries of the square}
Let $\theta\colon\kk\bS_N\to\kk\bS_N$ be the linear map $u\mapsto\sgn(u)\,u$, the \emph{sign twist}. Since $w_{L\leftarrow\sigma}$ differs from $w_{L\leftarrow\id}$ by a permutation of the positions in $L$ of sign $\sgn(\sigma)$, we have $\theta(g[w,L])=\pm(\theta g)[w,L]$ for every $g\in\kk\bS_k$. Hence $\theta$ maps the ideal generated by a set $G$ onto the ideal generated by $\theta(G)$. Moreover, $\theta$ is the automorphism of the algebra $\kk\bS_N$ given by tensoring with the sign representation. So it maps $\End(V_\lambda)$ onto $\End(V_{\lambda'})$ and $\cJ_\Lambda$ onto $\cJ_{\Lambda'}$, where $\Lambda'$ consists of the conjugate partitions. For example, $\theta(e_n)=s_n$, and $\theta$ maps $(e_n)$ onto $(s_n)$.

The dihedral group $D_4$ of order eight, the group of symmetries of the square, acts on $\bigsqcup_n\bS_n$; it is generated by the reverse $w\mapsto w^r$, the complement $w\mapsto w^c$ and the inverse $w\mapsto w^{-1}$. It preserves pattern containment, and each of its elements maps substitution maps to substitution maps, for instance $(w_{L\leftarrow\sigma})^{-1}=(w^{-1})_{w(L)\leftarrow\sigma^{-1}}$. On $\kk\bS_N$ the reverse and the complement are the right and the left multiplication by $w_0=N\cdots21$, and the inverse is the antiautomorphism $u\mapsto u^{-1}$, which maps every block to itself, since the representations $V_\lambda$ are self-dual. Hence every symmetric ideal is stable under $D_4$.

\subsubsection{Quotients}

\begin{proposition}\label{prop:rep}
For $n\ge1$ and $N\ge0$ there is an isomorphism of algebras and of $\bS_N\times\bS_N$-modules
\[
\bigl(\cS/(e_n)\bigr)(N)\cong\bigoplus_{\substack{\lambda\vdash N\\ \ell(\lambda)<n}}\End(V_\lambda),\qquad\text{so that}\qquad
\dim\bigl(\cS/(e_n)\bigr)(N)=\sum_{\substack{\lambda\vdash N\\ \ell(\lambda)<n}}(f^\lambda)^2 .
\]
\end{proposition}

\begin{proof}
The element $e_n$ spans the block of the sign representation $V_{(1^n)}$, and a partition contains $(1^n)$ if and only if it has at least $n$ rows. Now apply Proposition~\ref{prop:symideals}.
\end{proof}

For instance, $\dim(\cS/(e_2))(N)=1$ for all $N$, since only $\lambda=(N)$ has fewer than two rows; indeed, modulo $e_2=12-21$ all permutations of $[N]$ are equal.

Similarly, the element $s_m$ spans the block of the trivial representation $V_{(m)}$, and a partition contains $(m)$ if and only if $\lambda_1\ge m$. So by Proposition~\ref{prop:symideals}, for $m\ge1$ and $N\ge0$,
\begin{equation}\label{eq:sm}
\bigl(\cS/(s_m)\bigr)(N)\cong\bigoplus_{\substack{\lambda\vdash N\\ \lambda_1<m}}\End(V_\lambda),\qquad
\dim\bigl(\cS/(s_m)\bigr)(N)=\sum_{\substack{\lambda\vdash N\\ \lambda_1<m}}(f^\lambda)^2 .
\end{equation}

\begin{proposition}\label{prop:ensm}
For $n,m\ge1$ and $N\ge0$,
\[
\bigl(\cS/(e_n,s_m)\bigr)(N)\cong\bigoplus_{\substack{\lambda\vdash N\\ \ell(\lambda)<n,\ \lambda_1<m}}\End(V_\lambda),\qquad
\dim\bigl(\cS/(e_n,s_m)\bigr)(N)=\sum_{\substack{\lambda\vdash N\\ \ell(\lambda)<n,\ \lambda_1<m}}(f^\lambda)^2 .
\]
In particular $(\cS/(e_n,s_m))(N)=0$ for $N>(n-1)(m-1)$, so the algebra $\cS/(e_n,s_m)$ is finite-dimensional.
\end{proposition}

\begin{proof}
By Proposition~\ref{prop:symideals}, $(e_n,s_m)$ is the symmetric ideal of the union of the upper sets of $(1^n)$ and of $(m)$. The remaining partitions are those contained in the rectangle with $n-1$ rows and $m-1$ columns, which have at most $(n-1)(m-1)$ cells.
\end{proof}

By the Robinson--Schensted correspondence \cite{Schensted}, $\dim(\cS/(e_n,s_m))(N)$ is the number of permutations of $[N]$ avoiding both $n\cdots21$ and $12\cdots m$. So the vanishing in Proposition~\ref{prop:ensm} is an algebraic counterpart of the Erd\H{o}s--Szekeres theorem \cite{ES}.

A partition is a \emph{hook} if it is of the form $(a,1^b)$, or equivalently if its diagram does not contain the square $(2,2)$. More generally, for $k,l\ge0$ let $H(k,l)$ be the set of partitions $\lambda$ with $\lambda_{k+1}\le l$, that is, of the partitions whose diagram lies in the union of the first $k$ rows and the first $l$ columns. A partition lies outside $H(k,l)$ if and only if it contains the rectangle
\[
\rho_{k,l}=\bigl((l+1)^{k+1}\bigr)
\]
with $k+1$ rows and $l+1$ columns. In particular $\lambda\in H(k,l)$ if and only if $\lambda'\in H(l,k)$. Let $\cJ_{k,l}$ be the symmetric ideal generated by the block $\End(V_{\rho_{k,l}})\subset\kk\bS_{(k+1)(l+1)}$, and let $\cA_{k,l}=\cS/\cJ_{k,l}$. By Proposition~\ref{prop:symideals},
\begin{equation}\label{eq:hook-dim}
\cA_{k,l}(N)\cong\bigoplus_{\substack{\lambda\vdash N\\ \lambda\in H(k,l)}}\End(V_\lambda),\qquad
\dim\cA_{k,l}(N)=\sum_{\substack{\lambda\vdash N\\ \lambda\in H(k,l)}}(f^\lambda)^2 .
\end{equation}
For $l=0$ the rectangle is the column $(1^{k+1})$ and $\cJ_{k,0}=(e_{k+1})$; for $k=0$ it is the row $(l+1)$ and $\cJ_{0,l}=(s_{l+1})$. For $k=l=1$ the rectangle is the square $(2,2)$, so $H(1,1)$ is the set of hooks and $\cA_{1,1}$ is the sum of the blocks of the hooks. The sign twist $\theta$ maps $\cJ_{k,l}$ onto $\cJ_{l,k}$.

\subsection{\texorpdfstring{$\bigoplus_n\End(V^{\otimes n})$}{End(V^n)} as a commutative bishuffle algebra}\label{sec:End}

Let $V$ be a vector space of finite dimension $d$. For $\pi\in\bS_N$ let $\Pi_\pi\in\End(V^{\otimes N})$ be the operator moving the $i$-th tensor factor to the place $\pi(i)$,
\[
\Pi_\pi(v_1\otimes\dots\otimes v_N)=v_{\pi^{-1}(1)}\otimes\dots\otimes v_{\pi^{-1}(N)},
\]
so that $\Pi_{\pi\rho}=\Pi_\pi\Pi_\rho$. Let $\cB_V$ be the sequence of spaces $\cB_V(N)=\End(V^{\otimes N})$ with the products
\begin{equation}\label{eq:BV}
X\star_{\sigma,\tau}Y=\Pi_\sigma\circ(X\otimes Y)\circ\Pi_\tau^{-1},\qquad X\in\cB_V(m),\quad Y\in\cB_V(n),
\end{equation}
for all pairs of $(m,n)$-shuffles $\sigma,\tau$. As a bi-collection, $\cB_V(I,J)=\operatorname{Hom}(V^{\otimes J},V^{\otimes I})$, and the product $\cB_V(I_1,J_1)\otimes\cB_V(I_2,J_2)\to\cB_V(I_1\sqcup I_2,J_1\sqcup J_2)$ is the tensor product of operators, where $V^{\otimes I_1}\otimes V^{\otimes I_2}$ is identified with $V^{\otimes(I_1\sqcup I_2)}$ by putting every tensor factor at the place given by its index.

Let $e_1,\dots,e_d$ be a basis of $V$. Then $\cB_V(N)$ has the basis of matrix units $E_{\alpha,\beta}$, indexed by the pairs of words $\alpha,\beta\in[d]^N$: the operator $E_{\alpha,\beta}$ maps $e_{\beta_1}\otimes\dots\otimes e_{\beta_N}$ to $e_{\alpha_1}\otimes\dots\otimes e_{\alpha_N}$ and the other basis vectors to zero. By \eqref{eq:BV},
\begin{equation}\label{eq:BV-words}
E_{\alpha,\beta}\star_{\sigma,\tau}E_{\gamma,\delta}=E_{\alpha',\beta'},\qquad\alpha'_{\sigma(i)}=(\alpha\gamma)_i,\quad\beta'_{\tau(i)}=(\beta\delta)_i,
\end{equation}
where $\alpha\gamma$ and $\beta\delta$ are concatenations.

\begin{proposition}\label{prop:BV}
\begin{enumerate}
\item The products \eqref{eq:BV} make $\cB_V$ into a commutative bishuffle algebra.
\item The group $GL(V)$ acts on every $\cB_V(N)$ by conjugation with $g^{\otimes N}$, the products are $GL(V)$-equivariant, and the spaces of invariants form a commutative bishuffle subalgebra
\[
\cB_V^{GL(V)}=\bigoplus_{N\ge0}\End_{GL(V)}(V^{\otimes N})\subset\cB_V .
\]
\item The map $\varphi_V\colon\cS\to\cB_V$, $u\mapsto\Pi_u$, is a homomorphism of bishuffle algebras with image in $\cB_V^{GL(V)}$. In terms of matrix units, for $u\in\bS_N$,
\[
\varphi_V(u)=\sum_{c\in[d]^N}E_{c\circ u^{-1},\,c}.
\]
\end{enumerate}
\end{proposition}

\begin{proof}
First, by the description of $\cB_V$ as a bi-collection, its product is the tensor product of operators, so it is associative. For the commutativity, let $X\in\cB_V(m)$ and $Y\in\cB_V(n)$, and let $c\in\bS_{m+n}$ be the permutation from the proof of Corollary~\ref{cor:S-commutative}, so that $\sigma^{op}=\sigma c$ and $\tau^{op}=\tau c$. The operator $\Pi_c$ moves the first $n$ tensor factors to the last $n$ places and the last $m$ factors to the first $m$ places, keeping their order, so $\Pi_c(Y\otimes X)\Pi_c^{-1}=X\otimes Y$. Hence
\[
Y\star_{\sigma^{op},\tau^{op}}X=\Pi_\sigma\Pi_c(Y\otimes X)\Pi_c^{-1}\Pi_\tau^{-1}=\Pi_\sigma(X\otimes Y)\Pi_\tau^{-1}=X\star_{\sigma,\tau}Y .
\]

Second, $g^{\otimes(m+n)}=g^{\otimes m}\otimes g^{\otimes n}$ commutes with the operators $\Pi_\pi$. So conjugation by $g^{\otimes(m+n)}$ maps $X\star_{\sigma,\tau}Y$ to the product along $(\sigma,\tau)$ of the conjugates of $X$ by $g^{\otimes m}$ and of $Y$ by $g^{\otimes n}$. Thus the products are equivariant, and a product of invariant operators is invariant.

Finally, $\Pi_{u\oplus v}=\Pi_u\otimes\Pi_v$, so by Proposition~\ref{prop:products}
\[
\Pi_{u\star_{\sigma,\tau}v}=\Pi_\sigma\,\Pi_{u\oplus v}\,\Pi_\tau^{-1}=\Pi_\sigma(\Pi_u\otimes\Pi_v)\Pi_\tau^{-1}=\Pi_u\star_{\sigma,\tau}\Pi_v ,
\]
and $\Pi_u$ commutes with $g^{\otimes N}$. The formula with matrix units is the identity $\Pi_u(e_{c_1}\otimes\dots\otimes e_{c_N})=e_{c_{u^{-1}(1)}}\otimes\dots\otimes e_{c_{u^{-1}(N)}}$.
\end{proof}

In other words, to compute $\varphi_V(u)$ we color the points of $u$ by $d$ colors in all possible ways and record the colors of the values and of the positions. For $d=1$ all permutations have the same image.

\begin{theorem}\label{thm:SW}
Let $\dim V=d$. The kernel of $\varphi_V$ is the ideal $(e_{d+1})$, and its image is $\cB_V^{GL(V)}$. Hence there is an isomorphism of bishuffle algebras
\[
\cS/(e_{d+1})\cong\cB_V^{GL(V)}=\bigoplus_{N\ge0}\End_{GL(V)}(V^{\otimes N}).
\]
In particular, the algebra of invariants is an $\bS$-symmetric commutative bishuffle algebra.
\end{theorem}

\begin{proof}
On $\kk\bS_N$ the map $\varphi_V$ is the action of $\bS_N$ on $V^{\otimes N}$ by permuting the tensor factors. By Schur--Weyl duality, its image is the commutant $\End_{GL(V)}(V^{\otimes N})$, and its kernel is the sum of the blocks $\End(V_\lambda)$ with $\ell(\lambda)>d$. By Proposition~\ref{prop:rep}, this sum is $(e_{d+1})(N)$.
\end{proof}

\subsubsection{Super vector spaces}
Let now $V=V_0\oplus V_1$ be a super vector space with $\dim V_0=k$ and $\dim V_1=l$. The group $\bS_N$ acts on $V^{\otimes N}$ by operators $\Pi_\pi$ permuting the tensor factors with the Koszul sign rule, and the formula \eqref{eq:BV} defines products on $\cB_V(N)=\End(V^{\otimes N})$, where the tensor product of operators also follows the Koszul rule, $(X\otimes Y)(a\otimes b)=(-1)^{|Y||a|}Xa\otimes Yb$. The argument of Proposition~\ref{prop:BV} gives
\[
X\star_{\sigma,\tau}Y=(-1)^{|X||Y|}\,Y\star_{\sigma^{op},\tau^{op}}X
\]
for homogeneous $X$ and $Y$. So $\cB_V$ is a commutative bishuffle algebra in the category of super vector spaces, and its even part is a commutative bishuffle algebra in the sense of Definition~\ref{def:commutative}. The operators $\Pi_u$ are even, and $\varphi_V\colon u\mapsto\Pi_u$ is again a homomorphism of bishuffle algebras. The action of $\bS_N$ on $V^{\otimes N}$ commutes with the action of the general linear Lie superalgebra $\mathfrak{gl}(V)=\mathfrak{gl}(k|l)$. We write $\End_{\mathfrak{gl}(V)}(V^{\otimes N})$ for the space of operators supercommuting with $\mathfrak{gl}(V)$.

\begin{theorem}\label{thm:super-SW}
Let $V$ be a super vector space of dimension $k|l$. The image of $\varphi_V$ in arity $N$ is $\End_{\mathfrak{gl}(V)}(V^{\otimes N})$, and the kernel of $\varphi_V$ is the ideal $\cJ_{k,l}$. Hence the spaces of $\mathfrak{gl}(V)$-invariants form a commutative bishuffle subalgebra of $\cB_V$, and
\[
\cA_{k,l}=\cS/\cJ_{k,l}\cong\bigoplus_{N\ge0}\End_{\mathfrak{gl}(V)}(V^{\otimes N})\cong\bigoplus_{\lambda\in H(k,l)}\End(V_\lambda),
\]
the last sum running over the partitions that do not contain the rectangle $\rho_{k,l}=((l+1)^{k+1})$.
\end{theorem}

\begin{proof}
By the theorem of Sergeev and of Berele and Regev \cite{Sergeev,BR}, the image of $\kk\bS_N$ in $\End(V^{\otimes N})$ is the commutant of $\mathfrak{gl}(V)$, and the irreducible representations of $\bS_N$ occurring in $V^{\otimes N}$ are the $V_\lambda$ with $\lambda\in H(k,l)$. So the kernel of $\varphi_V$ in arity $N$ is the sum of the blocks $\End(V_\lambda)$ with $\lambda\notin H(k,l)$, which is $\cJ_{k,l}(N)$ by \eqref{eq:hook-dim}.
\end{proof}

For $l=0$ this is Theorem~\ref{thm:SW}. In particular $\dim\cA_{k,l}(N)\le(k+l)^{2N}$.

\section{Noetherianity questions}\label{sec:noeth}

We say that a commutative bishuffle algebra $\cA$ is \emph{Noetherian} if every ideal of $\cA$ is finitely generated; equivalently, every ascending chain of ideals of $\cA$ stabilizes. The ideals of a quotient $\cS/\cI$ are the images of the ideals of $\cS$ containing $\cI$, so a quotient of a Noetherian algebra is Noetherian. In this subsection we show that the symmetric ideals of $\cS$ are always finitely generated, that $\cS$ itself is not Noetherian, that $\cB_V$ and the quotients $\cS/(e_k)$ are Noetherian, and we conjecture that all $\bS$-symmetric quotients other than $\cS$ itself are Noetherian.

\subsection{Infinite antichains}

A set of permutations closed under taking patterns is a \emph{class}. A class is \emph{partially well-ordered} if it contains no infinite antichain for pattern containment; since containment is well-founded, this is the same as being well-quasi-ordered. It is classical that the set of all permutations is not partially well-ordered \cite{SB,AMR,Vatter}. The following explicit antichain is convenient for us.

\begin{example}\label{ex:antichain}
For $j\ge2$ let
\[
a_j=2,\,3,\,5,\,1,\ 7,\,4,\ 9,\,6,\ \dots,\ 2j+1,\,2j-2,\ \ 2j+2,\,2j+3,\,2j\ \in\bS_{2j+3},
\]
where the pairs $2i+1,2i-2$ appear for $i=3,\dots,j$. For instance $a_2=2351674$, $a_3=235174896$ and $a_4=2\,3\,5\,1\,7\,4\,9\,6\,10\,11\,8$. The permutations $a_j$ form an infinite antichain, and they all avoid $321$.

Indeed, let $G(u)$ be the inversion graph of a permutation $u$: its vertices are the values of $u$, and two values are adjacent if they form an inversion. If $u$ is a pattern of $v$, then $G(u)$ is isomorphic to an induced subgraph of $G(v)$. The graph $G(a_j)$ is the path $1,5,4,7,6,9,8,\dots,2j+1,2j$ together with the two vertices $2,3$ joined to $1$ and the two vertices $2j+2,2j+3$ joined to $2j$. Its only vertices of degree three are $1$ and $2j$. Suppose that $a_i$ is a pattern of $a_j$ with $i<j$, and let $H$ be the corresponding induced subgraph of $G(a_j)$. It has two vertices of degree three, which must be $1$ and $2j$, together with all their neighbours. Since $H$ is connected, it contains the path from $1$ to $2j$, so $H=G(a_j)$, which contradicts $2i+3<2j+3$. Finally, an occurrence of $321$ is a triangle in the inversion graph, and $G(a_j)$ has no triangles.
\end{example}

\begin{proposition}\label{prop:antichain}
If $A$ is an infinite antichain of permutations, then the ideal of $\cS$ generated by $A$ is not finitely generated. In particular, $\cS$ is not Noetherian.
\end{proposition}

\begin{proof}
Let $A=\{a_1,a_2,\dots\}$. By Corollary~\ref{cor:classical}, a permutation lies in a monomial ideal if and only if it contains one of the generators. Every element of $(A)$ is a combination of permutations containing elements of $A$, so finitely many elements of $(A)$ lie in $(a_1,\dots,a_M)$ for some $M$. But $a_{M+1}\notin(a_1,\dots,a_M)$, since it contains none of $a_1,\dots,a_M$. So $(A)$ is not finitely generated, and the chain $(a_1)\subset(a_1,a_2)\subset\cdots$ does not stabilize.
\end{proof}

\subsection{The bishuffle algebras $\cB_V$ and $\cS/(e_k)$ are Noetherian}
On the other hand, the quotients $\cS/(e_k)$ are Noetherian. The idea is to replace the permutations by a larger, better ordered set of monomials: the matrix units of the algebra $\cB_V$ of Section~\ref{sec:End}, that is, the pairs of words in a finite alphabet.

\begin{lemma}\label{lem:BV-noeth}
For every finite-dimensional vector space $V$, the bishuffle algebra $\cB_V$ is Noetherian.
\end{lemma}

\begin{proof}
First, we call the matrix units $E_{\alpha,\beta}$ \emph{monomials}; by \eqref{eq:BV-words}, a product of monomials is a monomial. Say that $E_{\alpha,\beta}$ divides $E_{\alpha'',\beta''}$ if $E_{\alpha'',\beta''}=E_{\alpha,\beta}\star_{\sigma,\tau}E_{\gamma,\delta}$ for some monomial $E_{\gamma,\delta}$ and shuffles $\sigma,\tau$. Since $\sigma$ and $\tau$ can be chosen independently, this holds if and only if $\alpha$ is a subword of $\alpha''$ and $\beta$ is a subword of $\beta''$, the two embeddings being arbitrary. By Higman's lemma \cite{Higman}, the words in the finite alphabet $[d]$ are well-quasi-ordered by the subword order, and hence so are the pairs of words, ordered componentwise. So every set of monomials closed under taking multiples has finitely many minimal elements.

Second, order the monomials of $\cB_V(N)$ lexicographically by the word $\alpha$ and then by the word $\beta$, and let $\lm(f)$ be the largest monomial of an element $f$. If $E_{\alpha,\beta}<E_{\alpha_1,\beta_1}$, then $E_{\alpha,\beta}\star_{\sigma,\tau}E_{\gamma,\delta}<E_{\alpha_1,\beta_1}\star_{\sigma,\tau}E_{\gamma,\delta}$: the first difference of the new words $\alpha'$, or, if $\alpha=\alpha_1$, of the new words $\beta'$, occurs at the image of the first difference, because $\sigma$ and $\tau$ are increasing on $[n]$. The same holds in the second argument. Hence $\lm(f\star_{\sigma,\tau}E)=\lm(f)\star_{\sigma,\tau}E$ for every element $f$ and monomial $E$.

Finally, let $\cJ$ be an ideal of $\cB_V$. By the second part, the leading monomials of the elements of $\cJ$ form a set closed under taking multiples, so by the first part it has finitely many minimal elements $\lm(f_1),\dots,\lm(f_r)$ with $f_i\in\cJ$. Every $0\ne f\in\cJ(N)$ has a leading monomial of the form $\lm(f_i)\star_{\sigma,\tau}E$, and subtracting a multiple of $f_i\star_{\sigma,\tau}E$ lowers the leading monomial. By induction over the finite set of monomials of $\cB_V(N)$, $\cJ=(f_1,\dots,f_r)$.
\end{proof}

\begin{theorem}\label{thm:noeth-ek}
For every $k\ge1$ the bishuffle algebras $\cS/(e_k)$ and $\cS/(s_k)$ are Noetherian.
\end{theorem}

\begin{proof}
First, for $k=1$ the ideal $(e_1)$ contains every permutation of positive length, and there is nothing to prove. Let $k=d+1\ge2$, let $V=\kk^d$, and let $\cA=\cB_V^{GL(V)}\cong\cS/(e_k)$ (Theorem~\ref{thm:SW}). Since $GL(V)$ is reductive and $\kk$ has characteristic zero, every $\cB_V(N)$ is a semisimple $GL(V)$-module, and the projection $R\colon \cB_V(N)\to\cA(N)$ onto the invariants commutes with all $GL(V)$-equivariant linear maps. For $a\in\cA(n)$ the map $Y\mapsto a\star_{\sigma,\tau}Y$ from $\cB_V(m)$ to $\cB_V(n+m)$ is equivariant, so
\[
R(a\star_{\sigma,\tau}Y)=a\star_{\sigma,\tau}R(Y).
\]

Second, let $\cJ$ be an ideal of $\cA$, and let $\cJ\cB_V$ be the ideal of $\cB_V$ generated by $\cJ$; since $\cB_V$ is commutative, it is spanned by the products $a\star_{\sigma,\tau}Y$ with $a\in\cJ$ and $Y\in \cB_V$. By Lemma~\ref{lem:BV-noeth}, $\cJ\cB_V$ is finitely generated, and since each of finitely many generators is a finite sum of such products, it is generated by finitely many elements $a_1,\dots,a_r\in\cJ$. Let $x\in\cJ$. Then $x=\sum_\nu a_{i_\nu}\star_{\sigma_\nu,\tau_\nu}Y_\nu$ with $Y_\nu\in \cB_V$, and applying $R$ gives
\[
x=R(x)=\sum_\nu a_{i_\nu}\star_{\sigma_\nu,\tau_\nu}R(Y_\nu)
\]
with $R(Y_\nu)\in\cA$. So $\cJ$ is generated by $a_1,\dots,a_r$, and $\cS/(e_k)$ is Noetherian.

Finally, the sign twist $\theta$ of Section~\ref{sec:symdef} maps the ideal generated by a set $G$ onto the ideal generated by $\theta(G)$, and $\theta(e_k)=s_k$. So it induces a bijection between the ideals containing $(e_k)$ and those containing $(s_k)$ that preserves finite generation.
\end{proof}

\begin{corollary}\label{cor:noeth-strip}
Let $\Lambda$ be an upper set of Young's lattice that contains a column $(1^k)$ or a row $(k)$. Then $\cS/\cJ_\Lambda$ is Noetherian.
\end{corollary}

\begin{proof}
The ideal $\cJ_\Lambda$ contains $(e_k)$, respectively $(s_k)$, so $\cS/\cJ_\Lambda$ is a quotient of a Noetherian algebra.
\end{proof}

\begin{remark}\label{rem:noeth-ek}
\begin{enumerate}
\item In $\cS/(e_3)$ the ideal generated by the images of the permutations $a_j$ of Example~\ref{ex:antichain} is finitely generated. So there is an $M$ with $a_j\in(e_3,a_2,\dots,a_M)$ for all $j$, although the $a_j$ are pairwise incomparable and avoid $321$. The proof gives no bound for $M$.
\item The proof follows the strategy of Sam and Snowden \cite{SS}: the bishuffle algebra $\cB_V$, whose monomials are well-quasi-ordered by Higman's lemma, plays the role of a Gr\"obner cover, and the passage to the invariants is Hilbert's classical argument with the Reynolds operator. In the language of Noetherianity up to symmetry \cite{HS}, $\cB_V(N)$ is the multilinear part of the polynomial ring in two families of variables $y_{ia}$, $z_{ja}$, with $i,j\in[N]$ and $a\in[d]$, the matrix unit $E_{\alpha,\beta}$ corresponding to $\prod_iy_{i\alpha_i}\prod_jz_{j\beta_j}$, and the monoid of increasing maps acts on the indices $i$ and $j$ separately. For $d=\infty$ Higman's lemma fails, in accordance with Proposition~\ref{prop:antichain}.
\item The same argument shows that finitely generated modules over $\cS/(e_k)$ are Noetherian: the monomials of a free $\cB_V$-module on a generator of arity $n$ are the pairs of words in the alphabet $[d]\sqcup\{\ast\}$ with exactly $n$ letters $\ast$ in each word, and these are again well-quasi-ordered by Higman's lemma. We do not develop modules over bishuffle algebras here.
\end{enumerate}
\end{remark}

\subsection{On Noetherianity of $\mathbb{S}$-symmetric ideals}

\begin{lemma}\label{lem:young-wqo}
	Young's lattice is well-quasi-ordered: in every infinite sequence $\lambda^{(1)},\lambda^{(2)},\dots$ of partitions there are $i<j$ with $\lambda^{(i)}\subseteq\lambda^{(j)}$. In particular, every upper set of Young's lattice has only finitely many minimal elements.
\end{lemma}

\begin{proof}
	Regard a partition as the word $\lambda_1\lambda_2\cdots\lambda_\ell$ of its parts in the alphabet $\mathbb N$, ordered as usual. By Higman's lemma \cite{Higman}, in every infinite sequence of such words there are $i<j$ and a strictly increasing map $\varphi$ with $\lambda^{(i)}_s\le\lambda^{(j)}_{\varphi(s)}$ for all $s$. Since $\varphi(s)\ge s$ and the parts of $\lambda^{(j)}$ weakly decrease, $\lambda^{(i)}_s\le\lambda^{(j)}_{\varphi(s)}\le\lambda^{(j)}_s$ for all $s$, that is, $\lambda^{(i)}\subseteq\lambda^{(j)}$. The minimal elements of an upper set form an antichain, which is finite by the first statement.
\end{proof}

\begin{proposition}\label{prop:sym-noeth}
	Every symmetric ideal of $\cS$ is finitely generated. More precisely, if $\Lambda$ is an upper set of Young's lattice with minimal elements $\mu^{(1)},\dots,\mu^{(r)}$, then $\cJ_\Lambda$ is generated by the blocks $\End(V_{\mu^{(1)}}),\dots,\End(V_{\mu^{(r)}})$. Consequently, every ascending chain of symmetric ideals of $\cS$ stabilizes, and in every $\bS$-symmetric commutative bishuffle algebra every symmetric ideal is finitely generated.
\end{proposition}

\begin{proof}
	First, by Lemma~\ref{lem:young-wqo} the upper set $\Lambda$ has finitely many minimal elements, and $\Lambda$ is the union of the upper sets $\{\lambda:\lambda\supseteq\mu^{(i)}\}$. By Proposition~\ref{prop:symideals}, the ideal generated by $\End(V_{\mu^{(i)}})$ is the symmetric ideal of this upper set, so $\cJ_\Lambda$ is generated by the $r$ blocks, which are finite-dimensional.
	
	Second, by Proposition~\ref{prop:symideals} an ascending chain of symmetric ideals is an ascending chain of upper sets $\Lambda_1\subseteq\Lambda_2\subseteq\cdots$. The union $\Lambda$ is an upper set; its finitely many minimal elements lie in some $\Lambda_m$, and then $\Lambda_m=\Lambda$.
	
	Finally, the symmetric ideals of $\cS/\cJ$ are the images of the symmetric ideals of $\cS$ containing $\cJ$.
\end{proof}

So symmetric quotients of $\cS$ satisfy the ascending chain condition for symmetric ideals. The question is what happens for arbitrary ideals, which need not be stable under the permutations of positions and of values.

\begin{lemma}\label{lem:sym-in-hook}
Every nonzero symmetric ideal $\cJ\subset\cS$ contains one of the ideals $\cJ_{k,l}$ of Section~\ref{sec:symdef}. Hence every $\bS$-symmetric commutative bishuffle algebra $\cS/\cJ$ with $\cJ\ne0$ is a quotient of one of the algebras $\cA_{k,l}$.
\end{lemma}

\begin{proof}
By Proposition~\ref{prop:symideals}, $\cJ=\cJ_\Lambda$ for a nonempty upper set $\Lambda$. Let $\mu\in\Lambda$ be a nonempty partition, and let $k+1=\ell(\mu)$ and $l+1=\mu_1$. Then $\mu$ is contained in the rectangle $\rho_{k,l}=((l+1)^{k+1})$, so $\rho_{k,l}\in\Lambda$ and $\cJ_{k,l}\subseteq\cJ$.
\end{proof}

\begin{conjecture}\label{conj:noeth-hooks}
For all $k,l\ge0$ the algebra $\cA_{k,l}=\cS/\cJ_{k,l}$ is Noetherian. Equivalently, by Lemma~\ref{lem:sym-in-hook}, every $\bS$-symmetric commutative bishuffle algebra other than $\cS$ itself is Noetherian.
\end{conjecture}

For $l=0$ and for $k=0$ the conjecture is Theorem~\ref{thm:noeth-ek}, since $\cJ_{k,0}=(e_{k+1})$ and $\cJ_{0,l}=(s_{l+1})$. The first open case is the algebra $\cA_{1,1}$, whose components are the sums of the blocks of the hooks. The proof of Theorem~\ref{thm:noeth-ek} does not apply directly. By Theorem~\ref{thm:super-SW}, $\cA_{k,l}$ embeds into the bishuffle algebra $\cB_V$ for the super vector space $V$ of dimension $k|l$. The proof of Lemma~\ref{lem:BV-noeth} applies to this algebra, whose monomials are, up to sign, the pairs of words in the alphabet $[k+l]$, and shows that it is Noetherian. But the general linear supergroup is not reductive, and there is no Reynolds operator to pass to the invariants. 

\begin{remark}\label{rem:specht}
Finite generation questions of this kind also appear in the theory of polynomial identities. The Specht problem asks whether every T-ideal of the free associative algebra, that is, every ideal stable under all endomorphisms, is finitely generated as a T-ideal \cite{Specht}. In characteristic zero the answer is positive by a theorem of Kemer \cite{Kemer87,Kemer91}.
\end{remark}

\section{Gr\"obner bases for commutative bishuffle algebras}\label{sec:GB}

\subsection{Monomial orderings}

\begin{definition}
A \emph{monomial ordering} on $\cS$ is a choice of a total order on each $\bS_n$ such that for all $w\in\bS_n$ and all $k$-subsets $L\subset[n]$ the substitution map
\[
\bS_k\to\bS_n,\qquad \sigma\mapsto w_{L\leftarrow\sigma},
\]
is strictly increasing.
\end{definition}

\begin{remark}\label{rem:compatible}
A monomial ordering is the same as a \emph{compatible ordering}, that is, an ordering of the monomials of $\cS$ that is compatible with the bishuffle multiplications: for all $(n,m)$-shuffles $\sigma,\tau$ and all permutations $u,u'\in\bS_n$ and $v,v'\in\bS_m$,
\[
u<u'\ \Longrightarrow\ u\star_{\sigma,\tau}v<u'\star_{\sigma,\tau}v,\qquad v<v'\ \Longrightarrow\ u\star_{\sigma,\tau}v<u\star_{\sigma,\tau}v' .
\]
Indeed, by Proposition~\ref{prop:products}, the map $u\mapsto u\star_{\sigma,\tau}v$ is the substitution map $u\mapsto w_{L\leftarrow u}$ for $L=\tau([n])$ and $w=\id_n\star_{\sigma,\tau}v$, and the map $v\mapsto u\star_{\sigma,\tau}v$ is the substitution map into the positions $\tau([n+1,n+m])$ of $w=u\star_{\sigma,\tau}\id_m$. Conversely, every substitution map is a product map, since $w_{L\leftarrow u}=u\star_{\pi_{w(L)},\pi_L}w|_{[n]\setminus L}$. This is the analogue of the compatibility of a monomial order on a polynomial ring with multiplication, and of the compatibility of monomial orderings of shuffle algebras with their products \cite{DK}. Since every $\bS_n$ is finite, no well-ordering condition is needed. Naturality is built in as well: an ordering of $\bS_n$ is transported to the monomials of $\cS(I,J)$, $|I|=|J|=n$, along the order-preserving bijections, whereas for the free algebra $F$ naturality has to be required (Section~\ref{sec:bishuffle}).
\end{remark}

The dihedral group $D_4$ of Section~\ref{sec:symdef} maps substitution maps to substitution maps. Hence if $<$ is a monomial ordering and $\gamma\in D_4$, then $u<^\gamma v\iff\gamma(u)<\gamma(v)$ is again a monomial ordering.

\begin{proposition}\label{prop:lex}
The lexicographic ordering $<_{\lex}$ of one-line words is a monomial ordering on $\cS$. Consequently so are its eight transforms $<_\lex^\gamma$, $\gamma\in D_4$. In particular the opposite of the lexicographic ordering, $<_{\lex}^c$, for which the leading monomial of an element is its lexicographically smallest permutation, is a monomial ordering.
\end{proposition}

\begin{proof}
If $\sigma<_\lex\sigma'$, then $w_{L\leftarrow\sigma}$ and $w_{L\leftarrow\sigma'}$ first differ at the position $l_k$, where $k$ is the first position at which $\sigma$ and $\sigma'$ differ, and the values there compare as $\sigma(k)<\sigma'(k)$.
\end{proof}

\begin{remark}\label{rem:monotone}
In every monomial ordering either $21>12$ or $12>21$. In the first case every permutation $u\ne w_0$ satisfies $u<u_{L\leftarrow 21}$ for an occurrence $L$ of $12$ in $u$, so $w_0=k\cdots21$ is the largest element of $\bS_k$ for all $k$; in the second case $\id$ is the largest element. Hence an element of $\kk\bS_k$ whose support is all of $\bS_k$ has a monotone leading monomial, and the direction of monotonicity does not depend on $k$.
\end{remark}

\subsection{Gr\"obner bases}

Fix a monomial ordering. For $0\ne f\in\kk\bS_n$ let $\lm(f)$ be its largest permutation and $\lc(f)$ the coefficient of $\lm(f)$. By the definition of a monomial ordering,
\begin{equation}\label{eq:lm}
\lm\bigl(g[w,L]\bigr)=w\qquad\text{whenever } L \text{ is an occurrence of } \lm(g) \text{ in } w .
\end{equation}

\begin{definition}
Let $\cI\subset\cS$ be an ideal. A subset $G\subset\cI$ is a \emph{Gr\"obner basis} of $\cI$ if for every $0\ne f\in\cI$ the permutation $\lm(f)$ contains $\lm(g)$ for some $g\in G$. A permutation is \emph{normal} if it avoids all $\lm(g)$, $g\in G$.
\end{definition}

If $L$ is an occurrence of $\lm(g)$ in a permutation $w$ appearing in $f$ with coefficient $c$, the \emph{reduction} $f\mapsto f-\frac{c}{\lc(g)}\,g[w,L]$ replaces $w$ by smaller permutations, by \eqref{eq:lm}. Hence normal permutations span $\cS/(G)$, and as usual we obtain:

\begin{proposition}\label{prop:GB-dim}
Let $G$ generate $\cI$. The following are equivalent: $G$ is a Gr\"obner basis of $\cI$; the normal permutations form a basis of $\cS/\cI$; $\dim(\cS/\cI)(n)=\#\Av_n(\lm(G))$ for all $n$. In this case
\[
f_{\cS/\cI}(t)=\sum_{n\ge0}\frac{\#\Av_n(\lm(G))}{(n!)^2}\,t^n .
\]
\end{proposition}

\subsection{The Buchberger criterion}

\begin{definition}
Let $g_1,g_2\in G$ with $p_i=\lm(g_i)\in\bS_{k_i}$. A \emph{small common multiple} of $p_1$ and $p_2$ is a triple $(w,L_1,L_2)$, where $w\in\bS_n$, $L_i$ is an occurrence of $p_i$ in $w$, $L_1\cup L_2=[n]$ and $L_1\cap L_2\ne\varnothing$; if $g_1=g_2$ we require $L_1\ne L_2$. Its \emph{S-polynomial} is
\[
S(w;L_1,L_2)=\frac{1}{\lc(g_1)}\,g_1[w,L_1]-\frac{1}{\lc(g_2)}\,g_2[w,L_2].
\]
\end{definition}

By \eqref{eq:lm}, all permutations of $S(w;L_1,L_2)$ are smaller than $w$. Since the two occurrences share a position, $n\le k_1+k_2-1$; in particular a finite set $G$ has only finitely many small common multiples.

\begin{theorem}\label{thm:buchberger}
A set $G$ is a Gr\"obner basis of $(G)$ if and only if every S-polynomial reduces to zero modulo $G$.
\end{theorem}

\begin{proof}
The S-polynomials lie in $(G)$, so the condition is necessary. Conversely, by the standard argument (Bergman's diamond lemma) it suffices to show the following. Let $L_1$ and $L_2$ be occurrences of $p_1=\lm(g_1)$ and $p_2=\lm(g_2)$ in a permutation $w$, and let $h=\frac1{\lc(g_1)}g_1[w,L_1]-\frac1{\lc(g_2)}g_2[w,L_2]$. Then $h=\sum_j c_j\,g_{j}[u_j,K_j]$ with $g_j\in G$, $K_j$ an occurrence of $\lm(g_j)$ in $u_j$, and $u_j<w$.

First, let $L_1\cap L_2=\varnothing$. Substitutions on disjoint sets of positions commute. Writing $g_1=\sum c_\sigma\sigma$ and $g_2=\sum d_\rho\rho$, we get
\[
\sum_{\sigma}c_\sigma\,g_2\bigl[w_{L_1\leftarrow\sigma},L_2\bigr]=\sum_{\sigma,\rho}c_\sigma d_\rho\,w_{L_1\leftarrow\sigma,\,L_2\leftarrow\rho}=\sum_\rho d_\rho\,g_1\bigl[w_{L_2\leftarrow\rho},L_1\bigr].
\]
Separating the terms $\sigma=p_1$ and $\rho=p_2$ gives
\[
\lc(g_1)\lc(g_2)\,h=\sum_{\sigma\ne p_1}c_\sigma\,g_2\bigl[w_{L_1\leftarrow\sigma},L_2\bigr]-\sum_{\rho\ne p_2}d_\rho\,g_1\bigl[w_{L_2\leftarrow\rho},L_1\bigr],
\]
and $w_{L_2\leftarrow\rho}<w$, $w_{L_1\leftarrow\sigma}<w$ by the definition of a monomial ordering.

Second, let $L_1\cap L_2\ne\varnothing$. Put $U=L_1\cup L_2$ and $u=w|_U$. The occurrences $L_1,L_2$ correspond to occurrences $L_1',L_2'$ in $u$, and $(u,L_1',L_2')$ is a small common multiple. By assumption $S(u;L_1',L_2')=\sum_j c_j\,g_j[u_j,K_j]$ with $u_j<u$. The map $x\mapsto w_{U\leftarrow x}$ is linear, strictly increasing, and transforms $g[x,K]$ into $g[w_{U\leftarrow x},K'']$, where $K''\subset U$ corresponds to $K$. Applying it to the last identity expresses $h$ in the required form.
\end{proof}

\begin{corollary}\label{cor:finite}
Let $G$ be a finite set of elements of arity at most $d$, generating an ideal $\cI$. Then $G$ is a Gr\"obner basis of $\cI$ if and only if $\dim(\cS/\cI)(n)=\#\Av_n(\lm(G))$ for all $n\le 2d-1$.
\end{corollary}

\begin{proof}
All S-polynomials lie in $\cI(n)$ with $n\le 2d-1$. If $\dim(\cS/\cI)(n)=\#\Av_n(\lm(G))$, then every nonzero element of $\cI(n)$ has a leading permutation containing some $\lm(g)$, so every element of $\cI(n)$ reduces to zero. Theorem~\ref{thm:buchberger} finishes the proof.
\end{proof}

Thus, for a fixed finite $G$, the Gr\"obner property is decidable by a finite computation. On the other hand, classical containment is not a well-quasi-order, so Buchberger's completion procedure need not terminate.

\subsection{Wilf equivalences from two orderings}

\begin{proposition}\label{prop:wilf}
If $G$ is a Gr\"obner basis of an ideal $\cI\subset\cS$ with respect to two monomial orderings $<$ and $<'$, then
\[
\#\Av_n\bigl(\lm_<(G)\bigr)=\dim(\cS/\cI)(n)=\#\Av_n\bigl(\lm_{<'}(G)\bigr)\qquad\text{for all } n .
\]
\end{proposition}

\begin{proof}
This is Proposition~\ref{prop:GB-dim} applied twice.
\end{proof}

For binomial ideals one obtains more. Let $p,q\in\bS_k$ and $f=p-q$. Let $\sim$ be the equivalence relation on $\bS_n$ generated by $w\sim w_{L\leftarrow q}$ for all occurrences $L$ of $p$ in $w$.

\begin{proposition}\label{prop:binomial}
\begin{enumerate}
\item $\dim(\cS/(f))(n)$ is the number of $\sim$-classes in $\bS_n$.
\item If $\lm_<(f)=p$ and $\{f\}$ is a Gr\"obner basis for $<$, then every $\sim$-class contains exactly one permutation avoiding $p$.
\item If moreover $\lm_{<'}(f)=q$ and $\{f\}$ is a Gr\"obner basis for $<'$, then sending a $p$-avoider to the unique $q$-avoider in its class is a bijection $\Av_n(p)\to\Av_n(q)$.
\end{enumerate}
\end{proposition}

\begin{proof}
The ideal $(f)(n)$ is spanned by the differences $f[w,L]=w_{L\leftarrow p}-w_{L\leftarrow q}$, which gives (1). In (2), replacing an occurrence of $p$ by $q$ decreases a permutation, so every class contains a $p$-avoider; by Proposition~\ref{prop:GB-dim} the number of classes equals $\#\Av_n(p)$, so each class contains exactly one. Part (3) follows from (2) applied to both orderings.
\end{proof}

\subsection{Examples: the quotients \texorpdfstring{$\cS/(e_n)$}{S/(e\_n)} and \texorpdfstring{$\cA_{1,1}$}{A\_(1,1)}}\label{sec:GB-examples}

We first record how the symmetries of Section~\ref{sec:symdef} interact with Gr\"obner bases. The sign twist $\theta$ multiplies every permutation by a sign, so it preserves leading monomials, and $G$ is a Gr\"obner basis if and only if $\theta(G)$ is. Let $\cJ$ be a symmetric ideal, $\gamma\in D_4$, and $<$ a monomial ordering. We have $\lm_{<^\gamma}(f)=\gamma^{-1}\bigl(\lm_<(\gamma f)\bigr)$, and $\cJ$ is stable under $\gamma$. So the leading permutations of $\cJ$ for $<^\gamma$ are the images under $\gamma^{-1}$ of its leading permutations for $<$. Since $\gamma$ preserves pattern containment, $\cJ$ has a Gr\"obner basis concentrated in arity $k$ for $<$ if and only if it has one for $<^\gamma$.

\subsubsection{The quotient $\cS/(e_n)$ and Schensted's theorem}

We use the Robinson--Schensted correspondence, a bijection between $\bS_N$ and the pairs of standard Young tableaux of the same shape $\lambda\vdash N$. By Schensted's theorem \cite{Schensted}, if $u$ corresponds to a pair of shape $\lambda$, then $\lambda_1$ is the length of a longest increasing subsequence of $u$, and $\ell(\lambda)$ is the length of a longest decreasing subsequence. Consequently
\begin{equation}\label{eq:schensted}
\#\Av_N(n\cdots21)=\sum_{\substack{\lambda\vdash N\\ \ell(\lambda)<n}}(f^\lambda)^2,\qquad
\#\Av_N(12\cdots n)=\sum_{\substack{\lambda\vdash N\\ \lambda_1<n}}(f^\lambda)^2,
\end{equation}
and the two numbers coincide, since $\lambda\mapsto\lambda'$ preserves $f^\lambda$.

\begin{theorem}\label{thm:en}
For every $n\ge1$, the antisymmetrizer $e_n$ is a Gr\"obner basis of the ideal $(e_n)$ for every monomial ordering of $\cS$. Its leading monomial is $n\cdots21$ if $21>12$, and $12\cdots n$ if $12>21$. In particular, $e_n$ is a Gr\"obner basis both for $<_\lex$, with leading monomial $n\cdots21$, and for $<_\lex^c$, with leading monomial $12\cdots n$.
\end{theorem}

\begin{proof}
First, the support of $e_n$ is all of $\bS_n$. By Remark~\ref{rem:monotone}, its leading monomial is the largest element of $\bS_n$, which is $n\cdots21$ if $21>12$ and $12\cdots n$ otherwise.

Second, in both cases $\#\Av_N(\lm(e_n))=\sum_{\ell(\lambda)<n}(f^\lambda)^2$ for all $N$, by Schensted's enumeration \eqref{eq:schensted}.

Finally, by Proposition~\ref{prop:rep} this number is the dimension of $(\cS/(e_n))(N)$. So $\dim(\cS/(e_n))(N)=\#\Av_N(\lm(e_n))$ for all $N$, and $e_n$ is a Gr\"obner basis by Proposition~\ref{prop:GB-dim}.
\end{proof}

\begin{remark}\label{rem:rsk}
By Proposition~\ref{prop:GB-dim}, the Gr\"obner property of $e_n$ for $<_\lex$ is equivalent to the enumeration \eqref{eq:schensted} of the permutations avoiding $n\cdots21$. So the argument can be reversed, and the Gr\"obner property gives a proof of Schensted's enumeration. Indeed, by Lemma~\ref{lem:minors}, $(e_n)(N)$ is the weight space of the ideal of $n\times n$ minors of the generic $N\times N$ matrix. These minors form a Gr\"obner basis for the antidiagonal and the diagonal term orders (Theorem~\ref{thm:ladder} for the full board), which gives the Gr\"obner property of $e_n$ for $<_\lex$ and $<_\lex^c$ by Lemma~\ref{lem:orders}. Sturmfels' proof of the Gr\"obner property of determinantal ideals \cite{Sturmfels} uses the Knuth correspondence \cite{Knuth}, whereas the proofs of Knutson, Miller and Yong \cite{KM,KMY} are geometric (Theorems~\ref{thm:KM} and~\ref{thm:KMY}). Through the latter, the enumeration \eqref{eq:schensted} has a proof that does not use the Robinson--Schensted correspondence. This determinantal argument is the case $\cW=\{\varnothing\}$ of Theorem~\ref{thm:bwx}.
\end{remark}

\begin{example}\label{ex:normal}
For $<_\lex$ the normal permutations of $\cS/(e_3)$ are the $321$-avoiders. A permutation containing $321$ is reduced by replacing an occurrence of $321$ by the alternating sum of the five other arrangements of its values; for instance
\[
321\equiv123-132-213+231+312\pmod{(e_3)} .
\]
The classes of the $321$-avoiders form a basis of $(\cS/(e_3))(N)$, whose dimension $\sum_{\ell(\lambda)\le2}(f^\lambda)^2$ is the Catalan number: $1,1,2,5,14,42,\dots$
\end{example}

\begin{remark}\label{rem:gessel}
The normalization by $(n!)^2$ built into bishuffle algebras is the one in which Gessel's theorem \cite{Gessel} is stated: the Hilbert series of $\cS/(e_n)$ is the Bessel determinant
\[
f_{\cS/(e_n)}(t)=\det\bigl[I_{|i-j|}(2\sqrt t)\bigr]_{1\le i,j\le n-1},\qquad I_j(2\sqrt t)=\sum_{m\ge0}\frac{t^{m+j/2}}{m!\,(m+j)!}.
\]
\end{remark}

\begin{problem}
Give a direct proof, by the Buchberger criterion in $\cS$, that $e_n$ is a Gr\"obner basis for all $n$.
\end{problem}

The same holds for the symmetrizer $s_m$.

\begin{proposition}\label{prop:sm}
For every $m\ge1$, the symmetrizer $s_m$ is a Gr\"obner basis of $(s_m)$ for every monomial ordering, with the same leading monomial as $e_m$.
\end{proposition}

\begin{proof}
The sign twist maps $e_m$ to $s_m$ and preserves leading monomials, so the statement follows from Theorem~\ref{thm:en}. Alternatively, the support of $s_m$ is all of $\bS_m$, so $\lm(s_m)=\lm(e_m)$ by Remark~\ref{rem:monotone}, and $\#\Av_N(\lm(s_m))=\sum_{\lambda_1<m}(f^\lambda)^2$ by Schensted's enumeration \eqref{eq:schensted}, which is the dimension \eqref{eq:sm} of the quotient.
\end{proof}

\begin{remark}\label{rem:ensm-notGB}
Let $n,m\ge2$. For every monomial ordering the leading monomials of $e_n$ and $s_m$ are monotone in the same direction, by Remark~\ref{rem:monotone}: they are $n\cdots21$ and $m\cdots21$ if $21>12$, and $12\cdots n$ and $12\cdots m$ if $12>21$. In the first case the identity permutations, and in the second case the decreasing permutations, avoid both leading monomials in every arity, whereas the quotient vanishes in large arity. So $\{e_n,s_m\}$ is a Gr\"obner basis of $(e_n,s_m)$ for no monomial ordering.
\end{remark}

\subsubsection{The algebra $\cA_{1,1}$}
Recall from Section~\ref{sec:symdef} that $\cA_{1,1}=\cS/\cJ_{1,1}$, where $\cJ_{1,1}$ is the symmetric ideal generated by the block $\End(V_{(2,2)})$, and that $\cA_{1,1}$ is the sum of the blocks of the hooks.

\begin{theorem}\label{thm:22}
Let $K=\{1234,2143,3412,4321\}\subset\bS_4$ be the Klein four-group, and for $u\in\bS_4$ let $[u]=\sum_{\kappa\in K}u\kappa$ be the sum of the elements of the coset $uK$.
\begin{enumerate}
\item The block $\End(V_{(2,2)})$ has the basis
\[
g_1=[4213]-[4132],\qquad g_2=[4231]-[4123],\qquad g_3=[4312]-[4123],\qquad g_4=[4321]-[4132];
\]
for instance $g_1=4213+2431+1342+3124-4132-1423-3241-2314$. So $\cJ_{1,1}$ is generated by $g_1,g_2,g_3,g_4$.
\item For $N\ge1$,
\[
\cA_{1,1}(N)\cong\bigoplus_{\substack{\lambda\vdash N\\ \lambda\text{ a hook}}}\End(V_\lambda),\qquad \dim\cA_{1,1}(N)=\binom{2N-2}{N-1}.
\]
\item The elements $g_1,g_2,g_3,g_4$ form a Gr\"obner basis of $\cJ_{1,1}$ for $<_\lex$, with the leading monomials $4213$, $4231$, $4312$ and $4321$. Consequently, the normal permutations are the permutations avoiding
\[
P_{1,1}=\{4213,\ 4231,\ 4312,\ 4321\},
\]
and $\#\Av_N(P_{1,1})=\binom{2N-2}{N-1}$ for all $N\ge1$.
\end{enumerate}
\end{theorem}

\begin{proof}
First, we describe the block. The group $K$ is normal in $\bS_4$ with quotient $\bS_3$, and the irreducible representations of $\bS_4$ on which $K$ acts trivially are $V_{(4)}$, $V_{(1^4)}$ and $V_{(2,2)}$, pulled back from the trivial, the sign and the two-dimensional representation of $\bS_3$. Hence the sum of these three blocks consists of the elements of $\kk\bS_4$ that are constant on the cosets $uK$, and $\End(V_{(2,2)})$ consists of those elements $f$ which moreover satisfy $\sum_uf_u=0=\sum_u\sgn(u)f_u$, where $f_u$ is the coefficient of $u$. Since $K$ consists of even permutations, the sign is constant on every coset, and in each $g_i$ the two cosets have the same sign. So $g_i\in\End(V_{(2,2)})$. Right multiplication by $K$ permutes the positions simply transitively, so every coset contains exactly one permutation beginning with $4$, and this permutation is its largest element for $<_\lex$. As $4\sigma>4\tau$ whenever $\sigma>\tau$, the leading monomials of $g_1,\dots,g_4$ are $4213$, $4231$, $4312$ and $4321$. They are distinct, so $g_1,\dots,g_4$ are linearly independent, and they form a basis of the four-dimensional block $\End(V_{(2,2)})$. This proves (1), since the ideal generated by a block is generated by any basis of it.

Second, by Proposition~\ref{prop:symideals} the quotient $\cA_{1,1}(N)$ is the sum of the blocks of the partitions of $N$ not containing $(2,2)$, which are the hooks $(N-b,1^b)$, $0\le b\le N-1$. A standard Young tableau of hook shape is determined by the set of entries in its first column other than $1$, so $f^{(N-b,1^b)}=\binom{N-1}{b}$, and $\sum_b\binom{N-1}{b}^2=\binom{2N-2}{N-1}$ by the Vandermonde identity. This proves (2).

Finally, the generators $g_1,\dots,g_4$ have arity four. By Corollary~\ref{cor:finite}, they form a Gr\"obner basis if and only if $\dim\cA_{1,1}(N)=\#\Av_N(P_{1,1})$ for all $N\le7$. For $N\le3$ both numbers equal $N!$: no permutation of length at most three contains a pattern of length four, and all partitions of $N\le3$ are hooks. For $N=4,5,6,7$ we count the normal permutations by a computer enumeration:
\[
\begin{array}{c|ccccccc}
N & 1 & 2 & 3 & 4 & 5 & 6 & 7\\\hline
\#\Av_N(P_{1,1}) & 1 & 2 & 6 & 20 & 70 & 252 & 924\\
\dim\cA_{1,1}(N)=\binom{2N-2}{N-1} & 1 & 2 & 6 & 20 & 70 & 252 & 924
\end{array}
\]
Hence $g_1,\dots,g_4$ form a Gr\"obner basis, and the rest of (3) follows from Proposition~\ref{prop:GB-dim}.
\end{proof}

\begin{remark}\label{rem:22-orderings}
Since $\cJ_{1,1}$ is stable under the symmetries of the square (Section~\ref{sec:symdef}), Theorem~\ref{thm:22} holds for each of the eight orderings $<_\lex^\gamma$, $\gamma\in D_4$: a basis of $\End(V_{(2,2)})$ with the four leading monomials $\gamma^{-1}(p)$, $p\in P_{1,1}$, is a Gr\"obner basis of $\cJ_{1,1}$ for $<_\lex^\gamma$. For instance, for $<_\lex^c$ the leading monomials are $1234$, $1243$, $1324$ and $1342$.
\end{remark}

\begin{remark}\label{rem:central-lit}
The complement maps $P_{1,1}$ onto the set of the patterns $1234$, $1243$, $1324$ and $1342$. The class of the permutations avoiding these four patterns is one of the twelve classes with a basis of four patterns of length four that Bruner lists as possibly enumerated by the central binomial coefficients \cite{Bruner}, and its generating function is recorded in the Permutation Pattern Avoidance Library \cite{PermPAL}, where it was obtained by an automated method. We have not found a written proof of the equality $\#\Av_N(P_{1,1})=\binom{2N-2}{N-1}$ in the literature. Theorem~\ref{thm:22} proves it and explains this enumeration: the permutations avoiding $P_{1,1}$ form a basis of $\cA_{1,1}(N)\cong\bigoplus\End(V_\lambda)$, the sum over the hooks $\lambda\vdash N$, and $\binom{2N-2}{N-1}$ is also the number of permutations of $[N]$ whose Robinson--Schensted shape is a hook. These two sets of permutations are different, and the second one is not even a permutation class: the permutation $25314$ has shape $(3,1,1)$, while its pattern $2413$ has shape $(2,2)$.
\end{remark}

\subsection{Gr\"obner bases and Noetherianity}\label{sec:GB-noeth}

We return to the Noetherianity questions of Section~\ref{sec:noeth}. A Gr\"obner basis whose normal permutations are well-quasi-ordered gives a direct proof of Noetherianity.

\begin{proposition}\label{prop:noeth-criterion}
\begin{enumerate}
\item Let $\cI\subset\cS$ be an ideal with a Gr\"obner basis for some monomial ordering, and suppose that its normal permutations form a partially well-ordered class. Then $\cS/\cI$ is Noetherian.
\item For a set $P$ of permutations, the monomial quotient $\cS/(P)$ is Noetherian if and only if $\Av(P)$ is partially well-ordered.
\end{enumerate}
\end{proposition}

\begin{proof}
First, let $C$ be the class of normal permutations of $\cI$, and let $\cJ\supseteq\cI$ be an ideal. Let $L$ be the set of leading permutations of the elements of $\cJ$ that lie in $C$. If $v=\lm(f)\in L$ with $f\in\cJ$, and $w\in C$ contains $v$ at the positions $K$, then $f[w,K]\in\cJ$ has the leading permutation $w$ by \eqref{eq:lm}, so $w\in L$. Hence $L$ is closed upwards in $C$, and since $C$ is partially well-ordered, $L$ has finitely many minimal elements $\lm(f_1),\dots,\lm(f_r)$ with $f_i\in\cJ$. We claim that $\cJ=\cI+(f_1,\dots,f_r)$. Let $f\in\cJ(N)$. Reducing $f$ modulo the Gr\"obner basis of $\cI$ we may assume that $f$ is a combination of normal permutations. If $f\ne0$, then $\lm(f)\in L$ contains some $\lm(f_i)$ at positions $K$, and subtracting a multiple of $f_i[\lm(f),K]$ replaces $f$ by an element of $\cJ$ with smaller leading permutation. Since $\bS_N$ is finite, the claim follows by induction. So $\cJ/\cI$ is generated by the images of $f_1,\dots,f_r$.

Finally, the set $P$ is a Gr\"obner basis of the monomial ideal $(P)$, with normal permutations $\Av(P)$, so (1) gives one direction of (2). Conversely, if $\Av(P)$ contains an infinite antichain $A$, then the argument of Proposition~\ref{prop:antichain}, applied to the monomial ideal $(P\cup A)$, shows that the ideal of $\cS/(P)$ generated by $A$ is not finitely generated.
\end{proof}

By the theorem of Atkinson, Murphy and Ru\v{s}kuc \cite{AMR}, the class $\Av(p)$ of a single permutation $p$ is partially well-ordered if and only if $p$ is one of $1$, $12$, $21$, $132$, $213$, $231$ and $312$. So $\cS/(p)$ is Noetherian exactly for these $p$; in particular $\cS/(321)$ and $\cS/(123)$ are not Noetherian, as Example~\ref{ex:antichain} shows directly.

Proposition~\ref{prop:noeth-criterion}(1) does not apply to the quotients $\cS/(e_k)$, $k\ge3$. For every monomial ordering, the leading monomial of $e_k$ is $k\cdots21$ or $12\cdots k$ (Remark~\ref{rem:monotone}), and by Theorem~\ref{thm:en} the normal permutations form the class $\Av(k\cdots21)$ or $\Av(12\cdots k)$. Both classes contain an infinite antichain: Example~\ref{ex:antichain} and its reverse. Nevertheless these quotients are Noetherian, by Theorem~\ref{thm:noeth-ek}. The same obstruction occurs for the algebra $\cA_{1,1}$ of Conjecture~\ref{conj:noeth-hooks}. Every element of $P_{1,1}=\{4213,4231,4312,4321\}$ contains $321$, so the class $\Av(P_{1,1})$ of normal permutations of $\cA_{1,1}$ contains the antichain of Example~\ref{ex:antichain}, and Proposition~\ref{prop:noeth-criterion}(1) does not apply; the same holds for the transforms of $<_\lex$ under the symmetries of the square (Remark~\ref{rem:22-orderings}).

\section{Matrix Schubert varieties and shape-Wilf equivalence of \texorpdfstring{$12\cdots k$ and $k\cdots21$}{12...k and k...21}}\label{sec:ferrers}

The space $\cS(n)=\kk\bS_n$ is the weight space of multidegree $(1,\dots,1;1,\dots,1)$ of the coordinate ring of $n\times n$ matrices, a permutation $u$ corresponding to the monomial $x_u=\prod_jx_{u(j),j}$. In Section~\ref{sec:examples} the relevant ideals were those of the determinantal varieties of the full matrix, which are stable under $GL_n\times GL_n$. In this section we impose determinantal conditions on Ferrers boards. For a Ferrers board $\lambda$, the ideal $I_k(\lambda)$ generated by the $k\times k$ minors with all entries in $\lambda$ is a ladder determinantal ideal, and its weight space is spanned by the transversals of $\lambda$. It is not an ideal of $\cS$, but ladder determinantal ideals are the building blocks of the ideals treated in Section~\ref{sec:BWX}.

Section~\ref{sec:boards} recalls Ferrers boards and pattern containment in a board. Section~\ref{sec:ladder} shows that the minors generating $I_k(\lambda)$ form Gr\"obner bases for a diagonal and for an antidiagonal term order. This is the only input from commutative algebra; it follows from the results of Knutson, Miller and Yong on matrix Schubert varieties. Section~\ref{sec:board-ideals} restates the result in terms of the antisymmetrizer $e_k$, which is a Gr\"obner basis on every Ferrers board, and deduces the shape-Wilf equivalence of $12\cdots k$ and $k\cdots21$.

\subsection{Ferrers boards and transversals}\label{sec:boards}

\subsubsection{Grids and Ferrers boards}
We draw grids like matrices: the cell $(i,j)$ lies in row $i$, counted from the top, and in column $j$, counted from the left. A permutation $u\in\bS_N$ is drawn as its permutation matrix $P_u$ of Section~\ref{sec:commutative}, the matrix of the linear map $e_j\mapsto e_{u(j)}$, and we identify $u$ with the set of its \emph{points} $(u(j),j)$, $j\in[N]$. Thus the columns are the positions and the rows are the values, and the pattern formed by a set of points is read from left to right. For cells $c=(a,b)$ and $c'=(a',b')$ we say that $c'$ is \emph{weakly north-west} of $c$, and write $c'\le c$, if $a'\le a$ and $b'\le b$; we say that $c'$ is \emph{strictly south-east} of $c$ if $a'>a$ and $b'>b$. The \emph{corner} of a finite set of cells is the cell $(a,b)$, where $a$ is the largest row and $b$ the largest column occurring in the set. It is the south-east corner of the bounding box of the set. The corner of an occurrence $L$ of a pattern in a permutation $u$ is the corner of the points $(u(l),l)$, $l\in L$.

\begin{definition}\label{def:ferrers}
Let $\rho$ and $\kappa$ be finite ordered sets, thought of as sets of rows and of columns. A subset $\lambda\subset\rho\times\kappa$ is a \emph{Ferrers board} if it is closed under moving north and west: if $(r,c)\in\lambda$, $r'\le r$ and $c'\le c$, then $(r',c')\in\lambda$. A cell of $\lambda$ is an \emph{outer corner} of $\lambda$ if neither the cell immediately below it nor the cell immediately to its right lies in $\lambda$.
\end{definition}

Let $\mu_r$ be the number of cells of $\lambda$ in row $r$. Then row $r$ of $\lambda$ consists of the first $\mu_r$ columns, and $\mu_r$ weakly decreases from top to bottom. Hence, after identifying $\rho$ and $\kappa$ with intervals, $\lambda$ is the Ferrers diagram, in English notation, of the partition formed by the nonzero row lengths \cite{Stanley}. The outer corners are the cells $(r,\mu_r)$ of the nonempty rows $r$ that are followed by a shorter row or are the last row, and $\lambda$ is the union of the rectangles $\{r'\le r\}\times\{c'\le c\}$ over its outer corners $(r,c)$. Ferrers boards are the classical setting of rook theory \cite{GJW}. Pattern avoidance on Ferrers boards is the main tool in the proofs of Wilf equivalences by Babson and West, by Stankova and West, and by Backelin, West and Xin \cite{BW,SW,BWX}.

\begin{definition}\label{def:transversal}
Let $\lambda\subset\rho\times\kappa$ be a Ferrers board with $|\rho|=|\kappa|$.
\begin{enumerate}
\item A \emph{transversal} of $\lambda$ is a bijection $T\colon\kappa\to\rho$ whose graph $\{(T(j),j):j\in\kappa\}$ lies in $\lambda$, that is, a placement of $|\rho|$ non-attacking rooks on $\lambda$.
\item Let $p\in\bS_k$. A transversal $T$ \emph{contains $p$ in $\lambda$} if some $k$ points of $T$ form the pattern $p$ and their corner lies in $\lambda$ \cite{BWX}. For $p=k\cdots21$ and $p=12\cdots k$ we say that $T$ has a \emph{decreasing}, respectively an \emph{increasing}, \emph{$k$-chain in $\lambda$}.
\item Two patterns are \emph{shape-Wilf-equivalent} if for every Ferrers board $\lambda$ they are avoided in $\lambda$ by the same number of transversals of $\lambda$ \cite{BW,BWX}.
\end{enumerate}
\end{definition}

Since $\lambda$ is closed under moving north and west, the corner of a set of points lies in $\lambda$ if and only if its whole bounding box does. For the full board $\lambda=\rho\times\kappa$, containment in $\lambda$ is ordinary pattern containment. The corner of an increasing chain is its last point, so a transversal has an increasing $k$-chain in $\lambda$ if and only if it contains $12\cdots k$ in the ordinary sense. For decreasing chains this is not the case.

\begin{figure}[ht]
\centering
\begin{tikzpicture}[scale=0.55]
\shaderows{1/5,2/5,3/3,4/3,5/2}
\gridlines{5}{5}
\draw[very thick] (0,0) -- (5,0) -- (5,-2) -- (3,-2) -- (3,-4) -- (2,-4) -- (2,-5) -- (0,-5) -- cycle;
\draw[thick,dashed] (3.08,-0.08) rectangle (4.92,-1.92);
\draw[thick,dotted] (1.08,-1.08) rectangle (3.92,-2.92);
\pt{5}{1}\pt{3}{2}\pt{4}{3}\pt{2}{4}\pt{1}{5}
\end{tikzpicture}
\caption{The Ferrers board with row lengths $5,5,3,3,2$ (shaded) and its transversal $53421$. The dashed box is the bounding box of the decreasing $2$-chain $(2,4),(1,5)$. The dotted box, spanned by $(3,2)$ and $(2,4)$, is not contained in the board.}\label{fig:board}
\end{figure}
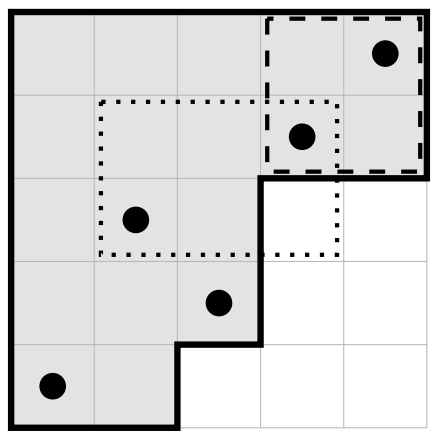

\begin{example}\label{ex:board}
Figure~\ref{fig:board} shows the Ferrers board $\lambda\subset[5]\times[5]$ with row lengths $5,5,3,3,2$, whose outer corners are $(2,5)$, $(4,3)$ and $(5,2)$, together with its transversal $T=53421$, which has the points $(5,1),(3,2),(4,3),(2,4),(1,5)$. The points $(2,4),(1,5)$ form a decreasing $2$-chain in $\lambda$, and so do $(5,1),(3,2)$. The points $(3,2),(2,4)$ form the pattern $21$ as well, but their corner $(3,4)$ is not in $\lambda$, so they do not form a chain in $\lambda$. The points $(3,2),(4,3)$ form an increasing $2$-chain in $\lambda$.
\end{example}

\subsubsection{Polynomial rings}
Let $\lambda\subset\rho\times\kappa$ be a Ferrers board and $R_\lambda=\kk[x_c:c\in\lambda]$, graded by $\mathbb Z^\rho\times\mathbb Z^\kappa$, where $x_{rc}$ has degree $(\varepsilon_r;\varepsilon_c)$. For the full board $[n]\times[n]$ we write $R_n$, and $Z=(x_{ij})$ for the generic matrix. For $A\subset\rho$ and $B\subset\kappa$ with $|A|=|B|$ let $\Delta_{A,B}$ be the minor of $Z$ with rows $A$ and columns $B$. Its corner is $(\max A,\max B)$, and $\Delta_{A,B}\in R_\lambda$ if and only if its corner lies in $\lambda$. Minors are homogeneous. For a bijection $T$ from a set of columns to a set of rows with graph in $\lambda$ we write $x_T=\prod_jx_{T(j),j}$, and for a set $X$ of cells we write $x_X=\prod_{c\in X}x_c$. If $|\rho|=|\kappa|$, the transversal monomials $x_T$ form a basis of the weight space
\[
\kk\lambda=(R_\lambda)_{(1;1)}
\]
of multidegree $(1,\dots,1;1,\dots,1)$. For $\lambda=[n]\times[n]$ we identify $\kk\lambda=(R_n)_{(1;1)}$ with $\kk\bS_n=\cS(n)$ via $x_u\leftrightarrow u$.

The substitutions of Section~\ref{sec:commutative} make sense for transversals. For a transversal $T$ of $\lambda$, a set $L$ of $k$ columns and $\sigma\in\bS_k$, let $T_{L\leftarrow\sigma}$ be the bijection that coincides with $T$ outside $L$ and maps $L$ onto $T(L)$ in the pattern $\sigma$, and for $g=\sum_\sigma c_\sigma\sigma\in\kk\bS_k$ put $g[T,L]=\sum_\sigma c_\sigma T_{L\leftarrow\sigma}$. If the corner of the points of $T$ in the columns $L$ lies in $\lambda$, then every $T_{L\leftarrow\sigma}$ is a transversal of $\lambda$, since its points in the columns $L$ lie in their bounding box.

\begin{lemma}\label{lem:minors}
Let $T$ be a transversal of a Ferrers board $\lambda$, let $L$ be a set of $k$ columns such that the corner of the points of $T$ in the columns $L$ lies in $\lambda$, and let $T'$ be the restriction of $T$ to the other columns. Then
\[
e_k[T,L]=\Delta_{T(L),L}\,x_{T'} .
\]
Conversely, every product $\Delta_{A,B}\,x_{T'}$ of a $k\times k$ minor with corner in $\lambda$ and a bijection $T'\colon\kappa\setminus B\to\rho\setminus A$ with graph in $\lambda$ is of this form.
\end{lemma}

\begin{proof}
Write $L=\{l_1<\dots<l_k\}$ and $T(L)=\{a_1<\dots<a_k\}$. The bijection $T_{L\leftarrow\sigma}$ maps $l_i$ to $a_{\sigma(i)}$, so $e_k[T,L]=\sum_\sigma\sgn(\sigma)\prod_ix_{a_{\sigma(i)},l_i}\,x_{T'}=\Delta_{T(L),L}\,x_{T'}$. Conversely, let $T$ consist of $T'$ and of the increasing bijection $B\to A$. Its graph lies in $\lambda$, because the corner of $\Delta_{A,B}$ does, and $e_k[T,B]=\Delta_{A,B}\,x_{T'}$.
\end{proof}

The weight space of a homogeneous ideal is spanned by the products of its generators with the monomials of complementary weight. Hence, for $\lambda=[n]\times[n]$, the component $(e_k)(n)$ of the ideal of $\cS$ generated by $e_k$ is the weight space of the ideal of $k\times k$ minors of the generic $n\times n$ matrix.

\subsection{Ladder determinantal ideals}\label{sec:ladder}

\begin{definition}\label{def:ladder}
Let $\lambda$ be a Ferrers board and $k\ge1$. The \emph{ladder determinantal ideal} $I_k(\lambda)\subset R_\lambda$ is the ideal generated by the $k\times k$ minors with corner in $\lambda$. The \emph{ladder determinantal variety} $Y_k(\lambda)$ is its zero set in the space $\Mat_\lambda$ of matrices indexed by $\lambda$.
\end{definition}

For $\lambda=[a]\times[b]$, $I_k(\lambda)$ is the ideal of $k\times k$ minors of the generic $a\times b$ matrix, and $Y_k(\lambda)$ is the variety of $a\times b$ matrices of rank less than $k$. In general, a $k\times k$ minor has its corner in $\lambda$ if and only if it is a minor of a north-west submatrix $Z_{[p]\times[q]}$ for an outer corner $(p,q)$ of $\lambda$. So $Y_k(\lambda)$ consists of the partial matrices $M$ whose north-west submatrices $M_{[p]\times[q]}$, for the outer corners $(p,q)$ of $\lambda$, have rank less than $k$. Ideals of this kind are the one-sided ladder determinantal ideals, see \cite[Section~2.4]{KM} and \cite{KMY} and the references there.

\subsubsection{Term orders}
A term order on a polynomial ring in entries of a matrix is \emph{antidiagonal} if the leading term of every minor is the product of its antidiagonal entries, and \emph{diagonal} if it is the product of its diagonal entries. Let $\lambda\subset\rho\times\kappa$ be a Ferrers board. Let $<_\ad$ be the lexicographic term order on $R_\lambda$ for the order of the variables by columns and, within a column, by decreasing rows, and let $<_\dg$ be the lexicographic term order for the order by columns and, within a column, by increasing rows. For $\lambda=[n]\times[n]$, the variables are ordered as $x_{n1}>x_{n-1,1}>\dots>x_{11}>x_{n2}>\dots>x_{1n}$, respectively as $x_{11}>x_{21}>\dots>x_{n1}>x_{12}>\dots>x_{nn}$.

\begin{lemma}\label{lem:orders}
The term order $<_\ad$ is antidiagonal and $<_\dg$ is diagonal. If $T$ and $T'$ are transversals of $\lambda$ and $j$ is the first column with $T(j)\ne T'(j)$, then $x_T>_\ad x_{T'}$ if and only if $T(j)>T'(j)$, and $x_T>_\dg x_{T'}$ if and only if $T(j)<T'(j)$. In particular, on $(R_n)_{(1;1)}=\kk\bS_n$ the term orders $<_\ad$ and $<_\dg$ restrict to the orderings $<_\lex$ and $<_\lex^c$ of Proposition~\ref{prop:lex}.
\end{lemma}

\begin{proof}
Consider a minor with rows $a_1<\dots<a_k$ and columns $b_1<\dots<b_k$. Each of its terms contains exactly one variable from the column $b_1$, and the largest of these for $<_\ad$ is $x_{a_kb_1}$. Among the terms containing $x_{a_kb_1}$, the largest variable from the column $b_2$ is $x_{a_{k-1}b_2}$, and so on; so the leading term is the antidiagonal one. Two transversal monomials $x_T$ and $x_{T'}$ have the same variables in the columns before $j$, and the larger of their variables in the column $j$ is $x_{\max(T(j),T'(j)),j}$; it divides the larger of the two monomials. For $<_\dg$ the rows are reversed.
\end{proof}

\subsubsection{Matrix Schubert varieties}
We use matrix Schubert varieties only to quote the Gr\"obner bases of ladder determinantal ideals. For $w\in\bS_n$ and $p,q\in[n]$ let $r_w(p,q)=\#\{j\le q:w(j)\le p\}$ be the number of points of $w$ in the rectangle $[p]\times[q]$. The \emph{matrix Schubert variety} of $w$ is
\[
X_w=\bigl\{M\in\Mat_n:\operatorname{rank}M_{[p]\times[q]}\le r_w(p,q)\text{ for all }p,q\bigr\},
\]
and $I_w\subset R_n$ is the ideal generated by the minors of size $r_w(p,q)+1$ of $Z_{[p]\times[q]}$, for all $p,q$ \cite{Fulton}. The \emph{Rothe diagram} of $w$ is $D(w)=\{(p,q):w^{-1}(p)>q\text{ and }w(q)>p\}$: draw from every point of $w$ a ray to the right and a ray downwards, then $D(w)$ consists of the cells met by no ray. Its \emph{essential set} $\Ess(w)$ consists of the cells $(p,q)\in D(w)$ such that neither $(p+1,q)$ nor $(p,q+1)$ lies in $D(w)$. The permutation $w$ is \emph{vexillary} if it avoids $2143$.

\begin{theorem}[Fulton \cite{Fulton}]\label{thm:fulton}
The ideal $I_w$ is generated by the minors of size $r_w(p,q)+1$ of $Z_{[p]\times[q]}$ for $(p,q)\in\Ess(w)$, the \emph{essential minors}.
\end{theorem}

\begin{theorem}[Knutson--Miller {\cite[Theorem~B]{KM}}]\label{thm:KM}
For every $w\in\bS_n$, the minors of size $r_w(p,q)+1$ of $Z_{[p]\times[q]}$, for all $p,q$, form a Gr\"obner basis of $I_w$ with respect to every antidiagonal term order.
\end{theorem}

\begin{theorem}[Knutson--Miller--Yong {\cite[Theorem~3.8]{KMY}}]\label{thm:KMY}
If $w$ is vexillary, then the essential minors form a Gr\"obner basis of $I_w$ with respect to every diagonal term order.
\end{theorem}

Depending on the convention for permutation matrices, $X_w$ is the matrix Schubert variety of $w$ or of $w^{-1}$ in \cite{KM,KMY}. Transposition maps one to the other and preserves the classes of diagonal and of antidiagonal term orders, so this choice does not matter here. Fulton identified the one-sided ladder determinantal ideals with the ideals $I_w$ of vexillary permutations \cite[Section~2.4]{KM}. For the ideals $I_k(\lambda)$ the following lemma makes this explicit.

\begin{lemma}\label{lem:board-vex}
Let $\lambda\subset[a]\times[b]$ be a Ferrers board with row lengths $\mu_1\ge\dots\ge\mu_a$, and let $k\ge1$. Let $\nu$ be the partition with parts $\nu_i=\mu_{i+k-1}-(k-1)$ for those $i\ge1$ for which this number is positive. Choose $M$ with $M\ge a$, $M\ge b$ and $M\ge k-1+\nu_i+i$ for all $i$. Let $v\in\bS_{M-k+1}$ be defined recursively by letting $v(i)$ be the smallest number larger than $\nu_i$ and different from $v(1),\dots,v(i-1)$, where $\nu_i=0$ for $i>\ell(\nu)$. Put $u=v^{-1}$ and $w=\id_{k-1}\oplus u\in\bS_M$, so that the point of $u$ in the row $i$ is $(i,v(i))$. Then:
\begin{enumerate}
\item $D(u)$ is the Ferrers diagram of $\nu$, and $D(w)$ is the set of cells of $\lambda$ lying in rows and columns at least $k$;
\item $\Ess(w)$ is the set of outer corners $(p,q)$ of $\lambda$ with $p,q\ge k$, and $r_w(c)=k-1$ for all $c\in D(w)$;
\item $w$ is vexillary, and $I_w=I_k(\lambda)\,R_M$ is generated by the $k\times k$ minors with corner in $\lambda$.
\end{enumerate}
\end{lemma}

\begin{proof}
First, the choice of $M$ guarantees that $v(i)$ exists. Let $j\le\nu_i$. Then $v(i)>\nu_i\ge j$, and $v(i')>\nu_{i'}\ge\nu_i\ge j$ for $i'<i$, so $u(j)>i$ and $(i,j)\in D(u)$. Let $j>\nu_i$. If $v(i)>j$, then $j$ is smaller than $v(i)$ and larger than $\nu_i$, so by the choice of $v(i)$ it equals some $v(i')$ with $i'<i$; hence $u(j)<i$. In both cases $(i,j)\notin D(u)$, so $D(u)$ is the diagram of $\nu$. The rays from the points $(i,i)$, $i<k$, of $w$ cover all cells in the first $k-1$ rows and columns, so $D(w)$ is the diagram of $u$ shifted by $(k-1,k-1)$. A cell $(p,q)$ with $p,q\ge k$ lies in $\lambda$ if and only if $q-(k-1)\le\nu_{p-k+1}$, which proves (1).

Second, a cell $(p,q)$ with $p,q\ge k$ has the same neighbours below and to the right in $D(w)$ as in $\lambda$, so the essential cells of $w$ are the outer corners of $\lambda$ in rows and columns at least $k$. For $c=(p,q)\in D(w)$, $r_w(c)$ is $k-1$ plus the number of points of $u$ in $[p-k+1]\times[q-k+1]$. This rectangle lies in the Ferrers board $D(u)$, which contains no point of $u$, so $r_w(c)=k-1$.

Finally, $u$ avoids $132$. Indeed, if $i<j<l$ and $u(i)<u(l)<u(j)$, then $(u(l),j)\in D(u)$, and the point $(u(i),i)$ lies weakly north-west of it, hence in the Ferrers board $D(u)$; but no point of $u$ lies in $D(u)$. In an occurrence of $2143$ in $w$ the first position is followed by a smaller value, so it is not among the first $k-1$ positions, and the occurrence lies in $u$. Since every occurrence of $2143$ contains one of $132$, $w$ is vexillary. By Theorem~\ref{thm:fulton} and (2), $I_w$ is generated by the $k\times k$ minors of the submatrices $Z_{[p]\times[q]}$ for the outer corners $(p,q)$ of $\lambda$ with $p,q\ge k$, and these are the $k\times k$ minors with corner in $\lambda$.
\end{proof}

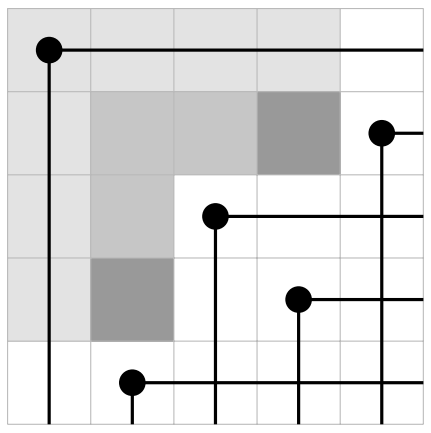
\begin{figure}[ht]
\centering
\begin{tikzpicture}[scale=0.55]
\shaderows{1/4,2/4,3/2,4/2}
\shadecell[gray!45]{2}{2}\shadecell[gray!45]{2}{3}\shadecell[gray!80]{2}{4}\shadecell[gray!45]{3}{2}\shadecell[gray!80]{4}{2}
\gridlines{5}{5}
\rays{1}{1}{5}\rays{2}{5}{5}\rays{3}{3}{5}\rays{4}{4}{5}\rays{5}{2}{5}
\pt{1}{1}\pt{2}{5}\pt{3}{3}\pt{4}{4}\pt{5}{2}
\end{tikzpicture}
\caption{The board $\lambda$ with row lengths $4,4,2,2$ (lightly shaded), $k=2$, and the vexillary permutation $w=15342$ of Lemma~\ref{lem:board-vex}. Its Rothe diagram consists of the cells of $\lambda$ in rows and columns at least two, and its essential cells (dark) are the outer corners $(2,4)$ and $(4,2)$ of $\lambda$, both of rank one.}\label{fig:vex}
\end{figure}

\begin{example}
For the board $\lambda$ with row lengths $4,4,2,2$ and $k=2$ one has $\nu=(3,1,1)$, $v=u=4231$ and $w=15342$ (Figure~\ref{fig:vex}). Thus the ideal generated by the $2\times2$ minors of the board is the ideal of $X_{15342}$, defined by $\operatorname{rank}M_{[2]\times[4]}\le1$ and $\operatorname{rank}M_{[4]\times[2]}\le1$.
\end{example}

\subsubsection{The Gr\"obner bases}

\begin{theorem}\label{thm:ladder}
Let $\lambda\subset\rho\times\kappa$ be a Ferrers board and $k\ge1$. The $k\times k$ minors with corner in $\lambda$ form a Gr\"obner basis of $I_k(\lambda)$ for every diagonal term order on $R_\lambda$ and for the antidiagonal term order $<_\ad$. Consequently, if $|\rho|=|\kappa|$, both
\[
\{x_T : T\text{ has no decreasing $k$-chain in }\lambda\}\quad\text{and}\quad\{x_T : T\text{ has no increasing $k$-chain in }\lambda\}
\]
are bases of $\bigl(R_\lambda/I_k(\lambda)\bigr)_{(1;1)}$. In particular, $12\cdots k$ and $k\cdots21$ are shape-Wilf-equivalent.
\end{theorem}

\begin{proof}
Identify $\rho$ and $\kappa$ with $[a]$ and $[b]$, let $w\in\bS_M$ be the vexillary permutation of Lemma~\ref{lem:board-vex}, and let $\mathcal G$ be the set of $k\times k$ minors with corner in $\lambda$. By Lemma~\ref{lem:board-vex}(3), $I_w=I_k(\lambda)\,R_M$ and $\mathcal G\subset I_w$. It suffices to show that $\mathcal G$ is a Gr\"obner basis of $I_w$ for a term order on $R_M$ that restricts to the given term order on $R_\lambda$. Indeed, then for $f\in I_k(\lambda)\subset I_w$ the leading term of $f$ is divisible by the leading term of an element of $\mathcal G$.

First, let $<$ be a diagonal term order on $R_\lambda$. We extend it to $R_M$ as follows: two monomials are compared first by their parts in the variables $x_c$, $c\notin\lambda$, using the lexicographic order for an order of these variables in which $x_{c'}<x_c$ whenever $c'\le c$ and $c'\ne c$, and then by $<$. This is a term order, and it is diagonal. Indeed, consider a minor with rows $a_1<\dots<a_s$ and columns $b_1<\dots<b_s$. If its corner lies in $\lambda$, the minor lies in $R_\lambda$, and its leading term is diagonal since $<$ is diagonal. Otherwise $x_{a_sb_s}$ is the largest of its variables outside $\lambda$, so its leading term is $x_{a_sb_s}$ times the leading term of the minor with rows $a_1,\dots,a_{s-1}$ and columns $b_1,\dots,b_{s-1}$, and the claim follows by induction on $s$. By Theorem~\ref{thm:KMY}, the essential minors form a Gr\"obner basis of $I_w$ for this term order. By Lemma~\ref{lem:board-vex}(2), they belong to $\mathcal G$, which is contained in $I_w$. Hence $\mathcal G$ is a Gr\"obner basis of $I_w$.

Second, consider $<_\ad$ on $R_M$; it restricts to $<_\ad$ on $R_\lambda$. By Theorem~\ref{thm:KM}, the initial ideal of $I_w$ is generated by the antidiagonal terms of minors $\Delta\in I_w$. Let $\Delta$ be such a minor and $C$ the set of cells of its antidiagonal, a decreasing chain. We claim that $C$ contains a decreasing $k$-chain whose corner lies in $\lambda$. Otherwise, let $P$ be the matrix with entries $1$ at the cells of $C$ and $0$ elsewhere. Then $\Delta(P)=\pm1$. A $k\times k$ minor with corner in $\lambda$ does not vanish at $P$ only if its rows and its columns are those of $k$ cells of $C$, and these cells would form a decreasing $k$-chain with corner in $\lambda$. So all elements of $\mathcal G$ vanish at $P$, hence so do all elements of $I_w=I_k(\lambda)R_M$, which contradicts $\Delta(P)\ne0$. Thus the antidiagonal term of $\Delta$ is divisible by the antidiagonal term of an element of $\mathcal G$, and $\mathcal G$ is a Gr\"obner basis of $I_w$.

Finally, $I_k(\lambda)$ is homogeneous for the row and column gradings, hence so is its initial ideal, and the leading monomials of the elements of $I_k(\lambda)_{(1;1)}$ are the transversal monomials lying in the initial ideal. A transversal monomial $x_T$ is divisible by the antidiagonal term of a minor with corner in $\lambda$ exactly when $T$ has a decreasing $k$-chain in $\lambda$, and by its diagonal term exactly when $T$ has an increasing $k$-chain in $\lambda$. Applied to $<_\ad$ and $<_\dg$, this gives the two bases.
\end{proof}

The last statement of Theorem~\ref{thm:ladder} is the core of the theorem of Backelin, West and Xin, who proved it by an explicit bijection \cite{BWX}; for $k=3$ it is due to Babson and West \cite{BW}. For the full board $\lambda=[a]\times[b]$, Theorem~\ref{thm:ladder} is the classical Gr\"obner basis theorem for determinantal ideals \cite{Sturmfels}.

\begin{remark}\label{rem:k2board}
Let $k=2$. By the remark after Definition~\ref{def:transversal}, a transversal has no increasing $2$-chain in $\lambda$ if and only if it is decreasing. The decreasing bijection $\kappa\to\rho$ sends the $i$-th column to the $i$-th row from the bottom. It is a transversal of $\lambda$ if and only if the $i$-th row from the bottom has at least $i$ cells for all $i$, and by Hall's condition this holds if and only if $\lambda$ has a transversal at all. Hence $\dim\bigl(R_\lambda/I_2(\lambda)\bigr)_{(1;1)}=1$ whenever $\lambda$ has a transversal. By Theorem~\ref{thm:ladder}, such a board also has exactly one transversal without a decreasing $2$-chain in $\lambda$.
\end{remark}

\begin{remark}\label{rem:projection}
The ladder determinantal variety is a projection, not an intersection. Let $\lambda\subset[a]\times[b]$ and let $I_k\subset R_{a\times b}$ be the ideal of $k\times k$ minors of the generic $a\times b$ matrix. Then
\[
I_k\cap R_\lambda=I_k(\lambda),
\]
so that, over an algebraically closed field, $Y_k(\lambda)$ is the closure of the image of the variety of $a\times b$ matrices of rank less than $k$ under the projection $\Mat_{a\times b}\to\Mat_\lambda$ that forgets the entries outside $\lambda$. Indeed, the lexicographic term order for an order of the variables in which every variable outside $\lambda$ is larger than every variable in $\lambda$, and $x_{c'}<x_c$ whenever $c'\le c$ and $c'\ne c$ lie on the same side, is diagonal by the argument in the proof of Theorem~\ref{thm:ladder}. So the $k\times k$ minors form a Gr\"obner basis of $I_k$ for it, by Theorem~\ref{thm:ladder} for the full board. It is an elimination order for the variables outside $\lambda$, so the minors lying in $R_\lambda$, which are those with corner in $\lambda$, form a Gr\"obner basis of $I_k\cap R_\lambda$. On the other hand, the intersection of the rank variety with the subspace of matrices vanishing outside $\lambda$ is in general smaller than $Y_k(\lambda)$. For the board $\lambda=\{(1,1),(1,2),(2,1)\}$ and $k=2$ one has $I_2(\lambda)=0$, while the matrices $\left(\begin{smallmatrix}x_{11}&x_{12}\\x_{21}&0\end{smallmatrix}\right)$ of rank less than two satisfy $x_{12}x_{21}=0$. This equation kills the only transversal monomial $x_{12}x_{21}$ of $\lambda$, which spans $\bigl(R_\lambda/I_2(\lambda)\bigr)_{(1;1)}$. It is the projection that matches the combinatorics of Babson and West, see Example~\ref{ex:tau1count}.
\end{remark}

\subsection{Board ideals and shape-Wilf equivalence}\label{sec:board-ideals}

Let $G$ be a set of elements of $\cS$, each of some arity. For a Ferrers board $\lambda\subset\rho\times\kappa$ with $|\rho|=|\kappa|$ let
\[
(G)_\lambda\subset\kk\lambda
\]
be the span of the elements $g[T,L]$, where $g\in G$ has arity $k$, $T$ is a transversal of $\lambda$, and $L$ is a set of $k$ columns such that the corner of the points of $T$ in the columns $L$ lies in $\lambda$. For the full board, $(G)_{[n]\times[n]}=(G)(n)$ is the component of the ideal generated by $G$. We order the transversals of $\lambda$ lexicographically, as in Lemma~\ref{lem:orders}: $T>_\lex T'$ if $T(j)>T'(j)$ for the first column $j$ with $T(j)\ne T'(j)$, and $<_\lex^c$ is the opposite ordering. As in \eqref{eq:lm}, $\lm(g[T,L])=T$ whenever the points of $T$ in the columns $L$ form an occurrence of $\lm(g)$ with corner in $\lambda$.

\begin{definition}\label{def:board-gb}
A set $G$ is a \emph{board Gr\"obner basis} for $<_\lex$ if for every Ferrers board $\lambda$ the leading transversals of the elements of $(G)_\lambda$ are exactly the transversals of $\lambda$ that contain $\lm(g)$ in $\lambda$ for some $g\in G$; similarly for $<_\lex^c$.
\end{definition}

By the remark before the definition, one inclusion always holds, so $G$ is a board Gr\"obner basis if and only if, for every $\lambda$, the dimension of $\kk\lambda/(G)_\lambda$ is the number of transversals of $\lambda$ avoiding all $\lm(g)$ in $\lambda$. Applied to the full boards, a board Gr\"obner basis is a Gr\"obner basis of the ideal $(G)$ of $\cS$. If $G$ is a board Gr\"obner basis for $<_\lex$ and for $<_\lex^c$, the two sets of leading monomials are shape-Wilf-equivalent.

\begin{corollary}\label{cor:shape-wilf}
For every $k\ge1$, the antisymmetrizer $e_k$ is a board Gr\"obner basis for $<_\lex$, with leading monomial $k\cdots21$, and for $<_\lex^c$, with leading monomial $12\cdots k$. Consequently, for every Ferrers board $\lambda$ the transversals of $\lambda$ avoiding $k\cdots21$ in $\lambda$ and the transversals of $\lambda$ avoiding $12\cdots k$ in $\lambda$ are two bases of the quotient $\kk\lambda/(e_k)_\lambda$, which for $\lambda=[n]\times[n]$ is the component $(\cS/(e_k))(n)$. In particular, $12\cdots k$ and $k\cdots21$ are shape-Wilf-equivalent.
\end{corollary}

\begin{proof}
By Lemma~\ref{lem:minors}, $(e_k)_\lambda=I_k(\lambda)_{(1;1)}$ is the weight space of the ladder determinantal ideal. So the statement follows from Theorem~\ref{thm:ladder} and Lemma~\ref{lem:orders}.
\end{proof}

\section{Backelin--West--Xin theorems from Gr\"obner bases}\label{sec:BWX}

Let $\cW$ be a nonempty set of permutations, possibly of different lengths, possibly infinite, and possibly containing the empty permutation $\varnothing$. We call its elements \emph{tails}. For $g=\sum_\sigma c_\sigma\sigma\in\kk\bS_k$ and a tail $\tau$ put $g\oplus\tau=\sum_\sigma c_\sigma\,\sigma\oplus\tau$, and write
\[
g\oplus\cW=\{g\oplus\tau:\tau\in\cW\}.
\]
In particular, for the antisymmetrizer and the symmetrizer,
\[
e_k\oplus\tau=\sum_{\sigma\in\bS_k}\sgn(\sigma)\,\sigma\oplus\tau,\qquad s_k\oplus\tau=\sum_{\sigma\in\bS_k}\sigma\oplus\tau;
\]
for instance, $e_2\oplus1=123-213$ and $e_k\oplus\varnothing=e_k$. In this section we compare the permutations avoiding the set of patterns $k\cdots21\oplus\cW$ with the permutations avoiding the set of patterns $12\cdots k\oplus\cW$. The set $\cW$ is fixed throughout the section, and all the constructions below depend on it. Let
\[
\cI_{k,\cW}=\bigl(e_k\oplus\cW\bigr)\subset\cS
\]
be the ideal generated by the elements $e_k\oplus\tau$, $\tau\in\cW$. Their leading monomials are $\lm_\lex(e_k\oplus\tau)=k\cdots21\oplus\tau$ and $\lm_{\lex^c}(e_k\oplus\tau)=12\cdots k\oplus\tau$.

\subsection{Overview}
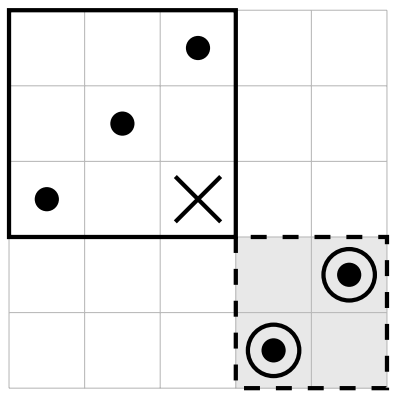
\begin{figure}[ht]
\centering
\begin{tikzpicture}[scale=0.5]
\fill[gray!18] (3,-3) rectangle (5,-5);
\gridlines{5}{5}
\draw[thick] (0,0) rectangle (3,-3);
\draw[thick,dashed] (3,-3) rectangle (5,-5);
\crosscell{3}{3}
\pt{1}{3}\pt{2}{2}\pt{3}{1}\pt{4}{5}\pt{5}{4}
\wit{4}{5}\wit{5}{4}
\end{tikzpicture}
\caption{The permutation $32154=321\oplus21$. The decreasing chain $(3,1),(2,2),(1,3)$ has the bounding box drawn in bold, with the corner $(3,3)$ (cross). The tail $21$ (circled) lies in the region strictly south-east of the corner (dashed), so it is a $\cW$-witness of the corner for every $\cW$ containing $21$.}\label{fig:occurrence}
\end{figure}

We follow the ideas of Babson and West \cite{BW} and of Backelin, West and Xin \cite{BWX}. An occurrence of $k\cdots21\oplus\tau$ in a permutation consists of a decreasing chain of $k$ points and of a copy of the tail $\tau$ lying strictly south-east of this chain (Figure~\ref{fig:occurrence}). So what matters is which cells of the grid have a copy of a tail from $\cW$ to the south-east. Starting from the set $\cW$, and independently of $k$, we introduce the following objects.
\begin{enumerate}
\item For $u\in\bS_N$, a copy of a tail from $\cW$ strictly south-east of a cell is a \emph{$\cW$-witness} of the cell. The cells having a $\cW$-witness form a Ferrers board $\lambda_\cW(u)$, the \emph{$\cW$-board} of $u$ (Section~\ref{sec:witness}).
\item The points of $u$ outside its $\cW$-board are its \emph{$\cW$-outer} points. The permutations with the same $\cW$-outer points form a \emph{$\cW$-class} $C$, and their remaining points run over all transversals of a Ferrers board $\lambda_C$ (Section~\ref{sec:classes}).
\end{enumerate}
The number $k$ enters only at the next step. The ideal $\cI_{k,\cW}$ splits along the $\cW$-classes, and on a class $C$ it becomes the weight space of the ladder determinantal ideal $I_k(\lambda_C)$ of Section~\ref{sec:ladder}. So, if $\Cl_\cW(N)$ denotes the set of $\cW$-classes in $\bS_N$, then by Theorem~\ref{thm:splitting}
\begin{equation}\label{eq:splitting}
\bigl(\cS/\cI_{k,\cW}\bigr)(N)\cong\bigoplus_{C\in\Cl_\cW(N)}\bigl(R_{\lambda_C}/I_k(\lambda_C)\bigr)_{(1;1)}.
\end{equation}
By Corollary~\ref{cor:shape-wilf}, each summand has two bases, formed by the transversals of $\lambda_C$ avoiding $k\cdots21$ in $\lambda_C$ and by those avoiding $12\cdots k$ in $\lambda_C$. Completed by the $\cW$-outer points of their classes, these transversals are exactly the permutations avoiding $k\cdots21\oplus\cW$, respectively $12\cdots k\oplus\cW$ (Section~\ref{sec:BWX-proof}). This leads to the main result of this section.

\begin{theorem}\label{thm:bwx}
Let $k\ge1$ and let $\cW$ be a nonempty set of permutations. Then $e_k\oplus\cW$ is a Gr\"obner basis of $\cI_{k,\cW}$ both for $<_\lex$ and for $<_\lex^c$, and so is $s_k\oplus\cW$ for the ideal it generates. Consequently, for all $N$,
\[
\#\Av_N\bigl(k\cdots21\oplus\cW\bigr)=\dim\bigl(\cS/\cI_{k,\cW}\bigr)(N)=\#\Av_N\bigl(12\cdots k\oplus\cW\bigr).
\]
\end{theorem}

For $\cW=\{\tau\}$ the equality of the numbers of avoiders is the theorem of Backelin, West and Xin \cite{BWX}, and for $k=3$ it is the theorem of Babson and West \cite{BW}. The extension to several tails is observed in \cite{BWX} right after Proposition~2.3, and it is stated as a theorem in \cite[Theorem~2]{BHKZ}.

\subsection{Witnesses and the board of a permutation}\label{sec:witness}

We use the pictures of Section~\ref{sec:boards}: a permutation $u\in\bS_N$ is the set of its points $(u(j),j)$ in the grid $[N]\times[N]$.

\begin{definition}\label{def:witness}
Let $u\in\bS_N$. A \emph{$\cW$-witness} of a cell $c\in[N]\times[N]$ in $u$ is a set of points of $u$, all strictly south-east of $c$, whose pattern lies in $\cW$. The \emph{$\cW$-board} of $u$ is the set of cells having a $\cW$-witness in $u$,
\[
\lambda_\cW(u)=\bigl\{c\in[N]\times[N]:c\text{ has a }\cW\text{-witness in }u\bigr\}.
\]
A point of $u$ is \emph{$\cW$-inner} if it lies in $\lambda_\cW(u)$, and \emph{$\cW$-outer} otherwise. We write $O_\cW(u)$ for the set of $\cW$-outer points of $u$.
\end{definition}

For $\cW=\{\tau\}$ the board $\lambda_\cW(u)$ is the board of Babson and West \cite{BW}. If $\varnothing\in\cW$, the empty set is a $\cW$-witness of every cell, so $\lambda_\cW(u)=[N]\times[N]$ and all points of $u$ are $\cW$-inner.

\begin{figure}[ht]
\centering
\begin{tikzpicture}[scale=0.5]
\fill[gray!18] (4,-2) rectangle (6,-6);
\gridlines{6}{6}
\draw[thick,dashed] (4,-2) rectangle (6,-6);
\crosscell{2}{4}
\pt{1}{4}\pt{2}{2}\pt{3}{6}\pt{4}{1}\pt{5}{5}\pt{6}{3}
\wit{3}{6}\wit{5}{5}
\end{tikzpicture}
\qquad
\begin{tikzpicture}[scale=0.5]
\fill[gray!18] (3,-4) rectangle (6,-6);
\gridlines{6}{6}
\draw[thick,dashed] (3,-4) rectangle (6,-6);
\crosscell{4}{3}
\pt{1}{4}\pt{2}{2}\pt{3}{6}\pt{4}{1}\pt{5}{5}\pt{6}{3}
\end{tikzpicture}
\qquad
\begin{tikzpicture}[scale=0.5]
\shaderows{1/4,2/4,3/2,4/2}
\gridlines{6}{6}
\pt{1}{4}\pt{2}{2}\pt{4}{1}
\opt{3}{6}\opt{5}{5}\opt{6}{3}
\end{tikzpicture}
\caption{$\cW$-witnesses for $\cW=\{21\}$ and $u=426153$. Left: the cell $(2,4)$ (cross), the region strictly south-east of it (dashed), and its $\cW$-witness $\{(3,6),(5,5)\}$ (circled). Middle: the cell $(4,3)$ has only one point strictly south-east of it, hence no $\cW$-witness. Right: the $\cW$-board $\lambda_\cW(u)$ (shaded), with the $\cW$-inner points black and the $\cW$-outer points white.}\label{fig:witness}
\end{figure}
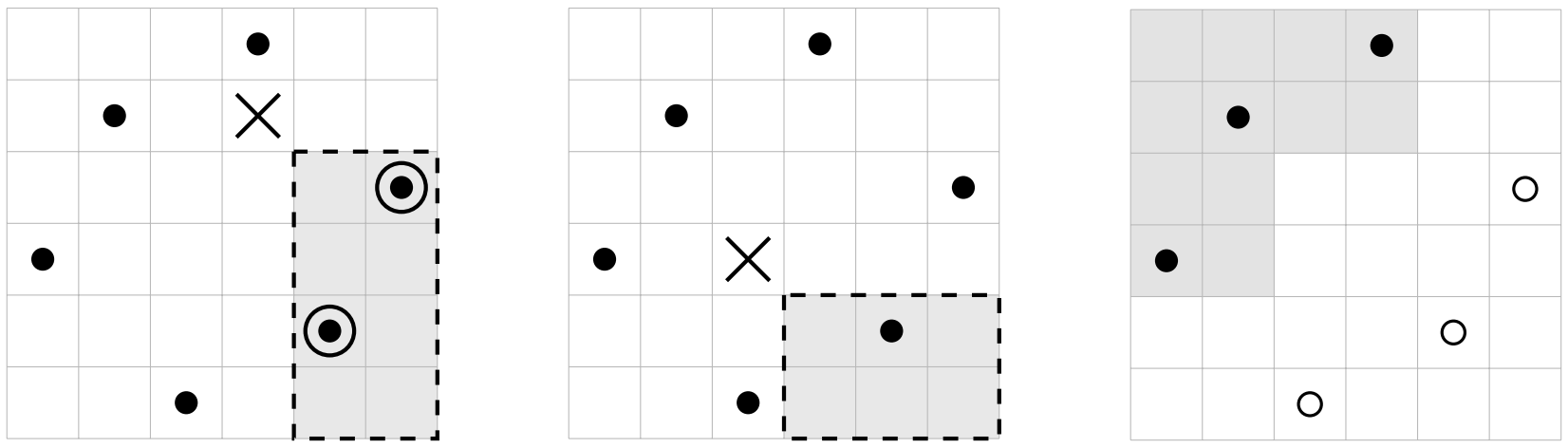

\begin{example}\label{ex:witness}
Let $\cW=\{21\}$ and $u=426153$ (Figure~\ref{fig:witness}). A $\cW$-witness of a cell $c$ is an inversion of $u$, that is, a pair of points forming the pattern $21$, lying strictly south-east of $c$.
\begin{itemize}
\item The cell $(2,4)$ has the $\cW$-witness $\{(3,6),(5,5)\}$, so $(2,4)\in\lambda_\cW(u)$. More generally, this set is a $\cW$-witness of every cell of $[2]\times[4]$, and $\{(5,5),(6,3)\}$ is a $\cW$-witness of every cell of $[4]\times[2]$.
\item The cell $(4,3)$ has the single point $(5,5)$ strictly south-east of it. So it has no $\cW$-witness, and $(4,3)\notin\lambda_\cW(u)$.
\end{itemize}
Hence $\lambda_\cW(u)=[2]\times[4]\cup[4]\times[2]$, with row lengths $4,4,2,2,0,0$. The $\cW$-inner points are $(1,4)$, $(2,2)$ and $(4,1)$, and the $\cW$-outer points are $(3,6)$, $(5,5)$ and $(6,3)$.
\end{example}

\begin{lemma}\label{lem:bw}
Let $u\in\bS_N$.
\begin{enumerate}
\item $\lambda_\cW(u)$ is a Ferrers board in $[N]\times[N]$.
\item Every cell of $\lambda_\cW(u)$ has a $\cW$-witness consisting of $\cW$-outer points of $u$. In particular, $\lambda_\cW(u)$ is determined by the set $O_\cW(u)$.
\item If the corner of a set of points of $u$ lies in $\lambda_\cW(u)$, then all these points are $\cW$-inner.
\end{enumerate}
\end{lemma}

\begin{proof}
First, if $c\in\lambda_\cW(u)$ and $c'\le c$, then every $\cW$-witness of $c$ is a $\cW$-witness of $c'$, so $c'\in\lambda_\cW(u)$.

Second, let $c\in\lambda_\cW(u)$. If $\varnothing\in\cW$, the empty set is a $\cW$-witness of $c$. Otherwise choose a $\cW$-witness $X$ of $c$ whose topmost row is as low as possible, and suppose that $X$ contains a $\cW$-inner point $x$. Then $x$ has a $\cW$-witness $X'$. Its points are strictly south-east of $x$, hence of $c$, so $X'$ is a $\cW$-witness of $c$ whose topmost row lies below the row of $x$, and hence below the topmost row of $X$. This contradicts the choice of $X$, so $X$ consists of $\cW$-outer points. Consequently, $\lambda_\cW(u)$ is the set of cells having a $\cW$-witness formed by points of $O_\cW(u)$.

Finally, every point of a set lies weakly north-west of its corner, so (3) follows from (1).
\end{proof}

\subsection{The \texorpdfstring{$\cW$}{W}-classes}\label{sec:classes}

\begin{definition}\label{def:class}
Two permutations $u,u'\in\bS_N$ are \emph{$\cW$-equivalent} if they have the same $\cW$-outer points, that is, if $O_\cW(u)=O_\cW(u')$. The equivalence classes are the \emph{$\cW$-classes}, and $\Cl_\cW(N)$ is the set of $\cW$-classes in $\bS_N$.
\end{definition}

Let $C\in\Cl_\cW(N)$, and let $O_C$ be the common set of $\cW$-outer points of its elements. By Lemma~\ref{lem:bw}(2), all elements of $C$ have the same $\cW$-board, which we denote by $\lambda_\cW(C)$. Let $\rho_C$ and $\kappa_C$ be the sets of rows and of columns not occupied by the points of $O_C$. The \emph{board of the class} $C$ is
\[
\lambda_C=\lambda_\cW(C)\cap(\rho_C\times\kappa_C),
\]
a Ferrers board in $\rho_C\times\kappa_C$ by Lemma~\ref{lem:bw}(1). For $u\in C$ let $T(u)$ be the set of $\cW$-inner points of $u$. These points lie in $\lambda_\cW(C)$ and occupy the rows $\rho_C$ and the columns $\kappa_C$, so $T(u)$ is a transversal of $\lambda_C$.

\begin{proposition}\label{prop:classes}
Let $C\in\Cl_\cW(N)$.
\begin{enumerate}
\item The map $u\mapsto T(u)$ is a bijection from $C$ onto the set of transversals of $\lambda_C$. Its inverse sends a transversal $T$ to the permutation with the points $O_C\sqcup T$.
\item For $u,u'\in C$ we have $u>_\lex u'$ if and only if $T(u)>_\lex T(u')$, in the ordering of transversals of Section~\ref{sec:board-ideals}.
\end{enumerate}
Consequently, for every $N$ there is a decomposition of vector spaces
\begin{equation}\label{eq:class-dec}
\kk\bS_N=\bigoplus_{C\in\Cl_\cW(N)}\kk C,\qquad\kk C\cong\kk\lambda_C,\quad u\mapsto T(u),
\end{equation}
which identifies the orderings $<_\lex$ and $<_\lex^c$ on $C$ with those on the transversals of $\lambda_C$.
\end{proposition}

\begin{proof}
First, the map is injective, since $u$ consists of the points $O_C$ and $T(u)$. Let $T$ be a transversal of $\lambda_C$, let $u'$ be the permutation with the points $O_C\sqcup T$, and let $u\in C$. By Lemma~\ref{lem:bw}(2), every cell of $\lambda_\cW(u)$ has a $\cW$-witness consisting of points of $O_C$, which are points of $u'$; so $\lambda_\cW(u)\subset\lambda_\cW(u')$. Conversely, let $c\notin\lambda_\cW(u)$, and suppose that $c$ has a $\cW$-witness $X$ in $u'$. If $X$ contained a point of $T$, then $c$ would lie north-west of a cell of $\lambda_\cW(u)$, hence in $\lambda_\cW(u)$ by Lemma~\ref{lem:bw}(1). So $X\subset O_C$, and $X$ is a $\cW$-witness of $c$ in $u$, a contradiction. Thus $\lambda_\cW(u')=\lambda_\cW(u)$. The points of $T$ lie in this board and the points of $O_C$ do not, so $O_\cW(u')=O_C$ and $T(u')=T$. This proves (1).

Second, two elements of $C$ coincide in the columns of the points of $O_C$, so the first column in which they differ belongs to $\kappa_C$. This proves (2).
\end{proof}

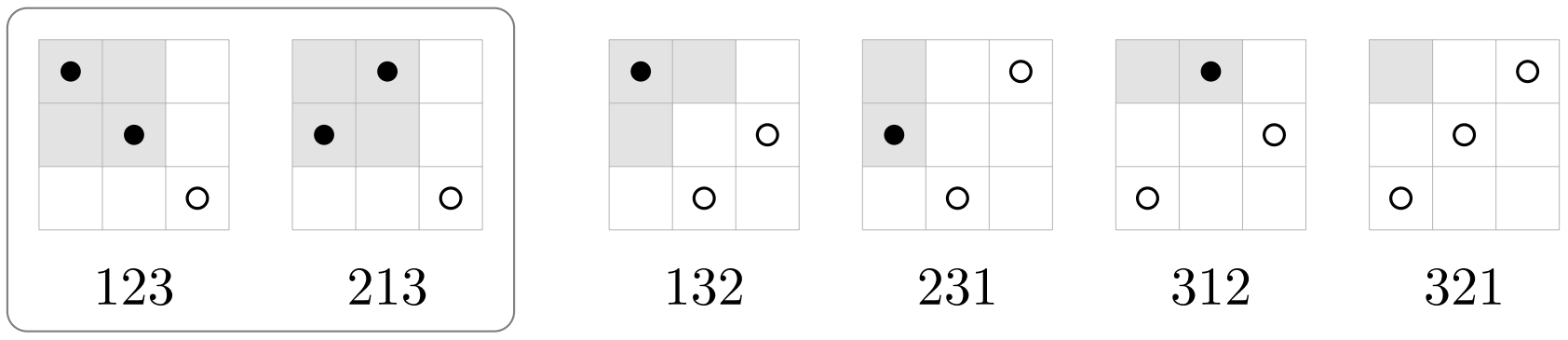
\begin{figure}[ht]
\centering
\begin{tikzpicture}[scale=0.45]
\draw[rounded corners,gray] (-0.5,0.5) rectangle (7.5,-4.6);
\begin{scope}
\shaderows{1/2,2/2}\gridlines{3}{3}
\pt{1}{1}\pt{2}{2}\opt{3}{3}
\node at (1.5,-3.9) {$123$};
\end{scope}
\begin{scope}[shift={(4,0)}]
\shaderows{1/2,2/2}\gridlines{3}{3}
\pt{1}{2}\pt{2}{1}\opt{3}{3}
\node at (1.5,-3.9) {$213$};
\end{scope}
\begin{scope}[shift={(9,0)}]
\shaderows{1/2,2/1}\gridlines{3}{3}
\pt{1}{1}\opt{2}{3}\opt{3}{2}
\node at (1.5,-3.9) {$132$};
\end{scope}
\begin{scope}[shift={(13,0)}]
\shaderows{1/1,2/1}\gridlines{3}{3}
\pt{2}{1}\opt{3}{2}\opt{1}{3}
\node at (1.5,-3.9) {$231$};
\end{scope}
\begin{scope}[shift={(17,0)}]
\shaderows{1/2}\gridlines{3}{3}
\pt{1}{2}\opt{2}{3}\opt{3}{1}
\node at (1.5,-3.9) {$312$};
\end{scope}
\begin{scope}[shift={(21,0)}]
\shaderows{1/1}\gridlines{3}{3}
\opt{1}{3}\opt{2}{2}\opt{3}{1}
\node at (1.5,-3.9) {$321$};
\end{scope}
\end{tikzpicture}
\caption{The five $\cW$-classes in $\bS_3$ for $\cW=\{1\}$. Each permutation is drawn with its $\cW$-board (shaded), the $\cW$-inner points black and the $\cW$-outer points white. The framed pair is the class $\{123,213\}$; the other classes have one element.}\label{fig:classes3}
\end{figure}

\begin{example}\label{ex:classes}
Let $\cW=\{1\}$. A cell lies in $\lambda_\cW(u)$ if and only if some point of $u$ lies strictly south-east of it, and a point of $u$ is $\cW$-outer if and only if it is a right-to-left maximum. So the $\cW$-class of $u$ consists of the permutations with the same right-to-left maxima, in the same positions.
\begin{enumerate}
\item For $N=3$ there are five classes, shown in Figure~\ref{fig:classes3}: the class $\{123,213\}$, whose board $\lambda_C$ is the full square on the rows and columns $\{1,2\}$, and four classes of size one. In $321$ all points are outer, although the cell $(1,1)$ lies in $\lambda_\cW(321)$.
\item For $N=4$ there are fourteen classes. Eight of them have size one, among them $\{3412\}$, whose board $\lambda_C\subset\{1,3\}\times\{1,3\}$ is the staircase $\{(1,1),(1,3),(3,1)\}$. The other six are listed in Table~\ref{tab:classes4}; for all of them $\lambda_C$ is a full square.
\item For $N=5$ there are $42$ classes, and boards that are not squares occur in classes with several elements. For instance, the class with the outer points $(3,5)$ and $(5,3)$ has the board $\lambda_C\subset\{1,2,4\}\times\{1,2,4\}$ with row lengths $3,3,2$, and it consists of $14523$, $24513$, $41523$ and $42513$.
\end{enumerate}
\end{example}

\begin{table}[ht]
\[
\begin{array}{c|c|c}
\text{outer points} & \rho_C\times\kappa_C & \text{class }C\\\hline
(4,4) & \{1,2,3\}\times\{1,2,3\} & 1234,\ 1324,\ 2134,\ 2314,\ 3124,\ 3214\\
(4,3),\ (3,4) & \{1,2\}\times\{1,2\} & 1243,\ 2143\\
(4,3),\ (2,4) & \{1,3\}\times\{1,2\} & 1342,\ 3142\\
(4,3),\ (1,4) & \{2,3\}\times\{1,2\} & 2341,\ 3241\\
(4,2),\ (3,4) & \{1,2\}\times\{1,3\} & 1423,\ 2413\\
(4,1),\ (3,4) & \{1,2\}\times\{2,3\} & 4123,\ 4213
\end{array}
\]
\caption{The $\cW$-classes with more than one element in $\bS_4$ for $\cW=\{1\}$. In each case the board $\lambda_C$ is all of $\rho_C\times\kappa_C$.}\label{tab:classes4}
\end{table}

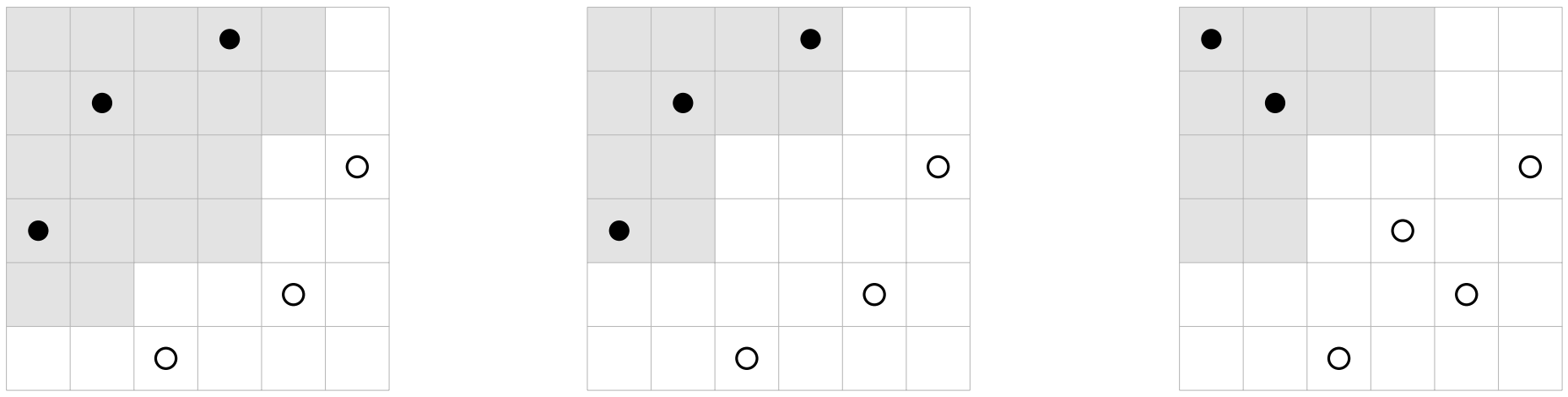
\begin{figure}[ht]
\centering
\begin{tikzpicture}[scale=0.5]
\shaderows{1/5,2/5,3/4,4/4,5/2}
\gridlines{6}{6}
\pt{1}{4}\pt{2}{2}\pt{4}{1}
\opt{3}{6}\opt{5}{5}\opt{6}{3}
\node at (3,-6.8) {\small $\cW=\{1\}$, $u=426153$};
\end{tikzpicture}
\qquad
\begin{tikzpicture}[scale=0.5]
\shaderows{1/4,2/4,3/2,4/2}
\gridlines{6}{6}
\pt{1}{4}\pt{2}{2}\pt{4}{1}
\opt{3}{6}\opt{5}{5}\opt{6}{3}
\node at (3,-6.8) {\small $\cW=\{21\}$, $u=426153$};
\end{tikzpicture}
\qquad
\begin{tikzpicture}[scale=0.5]
\shaderows{1/4,2/4,3/2,4/2}
\gridlines{6}{6}
\pt{1}{1}\pt{2}{2}
\opt{3}{6}\opt{4}{4}\opt{5}{5}\opt{6}{3}
\node at (3,-6.8) {\small $\cW=\{21\}$, $u=126453$};
\end{tikzpicture}
\caption{The dependence on $\cW$ in Example~\ref{ex:W-dependence}. Each permutation is drawn with its $\cW$-board (shaded), the $\cW$-inner points black and the $\cW$-outer points white. For $\cW=\{1\}$ the permutations $426153$ and $126453$ lie in one class; for $\cW=\{21\}$ they lie in different classes.}\label{fig:W-dependence}
\end{figure}

\begin{example}\label{ex:W-dependence}
The classes depend on $\cW$. Let $u=426153$ be the permutation of Example~\ref{ex:witness}. For both $\cW=\{1\}$ and $\cW=\{21\}$ its outer points are $(3,6)$, $(5,5)$ and $(6,3)$, but its boards are different (Figure~\ref{fig:W-dependence}).
\begin{itemize}
\item For $\cW=\{1\}$ the $\cW$-board of $u$ has row lengths $5,5,4,4,2,0$. The board of the class is the full square $\{1,2,4\}\times\{1,2,4\}$, and the class consists of the six permutations $126453$, $146253$, $216453$, $246153$, $416253$ and $426153$.
\item For $\cW=\{21\}$ the $\cW$-board of $u$ has row lengths $4,4,2,2,0,0$. The cell $(4,4)$ is no longer in it, so the board of the class has row lengths $3,3,2$, and the class consists of $146253$, $246153$, $416253$ and $426153$. In the two remaining permutations $126453$ and $216453$ the point $(4,4)$ has no $\cW$-witness. So it is $\cW$-outer, and these two permutations form another $\cW$-class, with the outer points $(3,6)$, $(4,4)$, $(5,5)$ and $(6,3)$.
\end{itemize}
Accordingly, the numbers of classes differ. For $N=3,4,5,6$ there are $5$, $14$, $42$, $132$ classes for $\cW=\{1\}$, and $6$, $23$, $103$, $513$ classes for $\cW=\{21\}$; for $\cW=\{\varnothing\}$ the whole of $\bS_N$ is one class. For $\cW=\{21\}$ and $N=4$ the only class with more than one element is $\{1243,2143\}$.
\end{example}

\begin{remark}\label{rem:young-classes}
Let $Y(u)$ be the minimal Young diagram containing all points of $u$, that is, the set of cells lying weakly north-west of some point of $u$. Its outer corners are the right-to-left maxima of $u$, so for $\cW=\{1\}$ the $\cW$-class of $u$ is the set of all permutations with the same diagram $Y(u)$. For other sets of tails the classes are in general smaller. If $\varnothing\notin\cW$, every $\cW$-class is contained in one class of equal minimal Young diagram. Indeed, a $\cW$-inner point has a nonempty $\cW$-witness, so it is not a right-to-left maximum. Hence the right-to-left maxima of $u$ are $\cW$-outer, and they are exactly those points of $O_\cW(u)$ with no point of $O_\cW(u)$ strictly south-east of them; so they are determined by $O_\cW(u)$. Example~\ref{ex:W-dependence} shows how a class of equal minimal Young diagram splits: some of its points become frozen, and the board shrinks.
\end{remark}

\subsection{The splitting of the ideal}\label{sec:splitting}

We now fix also $k\ge1$, and we describe the ideal $\cI_{k,\cW}$ by means of $\cW$-witnesses.

\begin{lemma}\label{lem:bwx-ideal}
For every $N$, the component $\cI_{k,\cW}(N)$ is spanned by the elements $e_k[w,A]$, where $w\in\bS_N$ and $A\subset[N]$ is a set of $k$ columns such that the corner of the points of $w$ in the columns $A$ has a $\cW$-witness in $w$.
\end{lemma}

\begin{proof}
First, let $\tau\in\cW$ have length $t$, let $w\in\bS_N$ and $L\subset[N]$ with $|L|=k+t$, and write $L=A\sqcup A'$, where $A$ consists of the $k$ smallest elements of $L$. Put $w'=w_{L\leftarrow\id_k\oplus\tau}$. For $\sigma\in\bS_k$, the permutation $w_{L\leftarrow\sigma\oplus\tau}$ coincides with $w'$ outside $A$ and takes on $A$ the $k$ smallest values of $w(L)$ in the pattern $\sigma$. So $w_{L\leftarrow\sigma\oplus\tau}=w'_{A\leftarrow\sigma}$, and
\[
(e_k\oplus\tau)[w,L]=e_k[w',A].
\]
The points of $w'$ in the columns $A'$ form the pattern $\tau$, and they lie to the right of the columns $A$ and below the rows $w'(A)$. So they form a $\cW$-witness of the corner of the points of $w'$ in the columns $A$.

Second, let $w$ and $A$ be as in the statement, let $X$ be a $\cW$-witness of the corner, with pattern $\tau\in\cW$, let $A'$ be the set of columns of $X$, and put $L=A\sqcup A'$. The points of $X$ lie to the right of the columns $A$ and below the rows $w(A)$. So $w_{A\leftarrow\sigma}=w_{L\leftarrow\sigma\oplus\tau}$ for all $\sigma\in\bS_k$, and $e_k[w,A]=(e_k\oplus\tau)[w,L]$.
\end{proof}

\begin{theorem}\label{thm:splitting}
For every $N$,
\[
\cI_{k,\cW}(N)=\bigoplus_{C\in\Cl_\cW(N)}\cI_C,\qquad\text{where}\quad\cI_C=\cI_{k,\cW}(N)\cap\kk C,
\]
and the isomorphism $\kk C\cong\kk\lambda_C$ of \eqref{eq:class-dec} maps $\cI_C$ onto the board ideal $(e_k)_{\lambda_C}=I_k(\lambda_C)_{(1;1)}$. In particular, the splitting \eqref{eq:splitting} holds.
\end{theorem}

\begin{proof}
First, let $e_k[w,A]$ be one of the elements of Lemma~\ref{lem:bwx-ideal}, let $c$ be the corner of the points of $w$ in the columns $A$, and let $C$ be the $\cW$-class of $w$. Then $c\in\lambda_\cW(C)$, and by Lemma~\ref{lem:bw}(3) the points of $w$ in the columns $A$ are $\cW$-inner. So $A\subset\kappa_C$ and $w(A)\subset\rho_C$, and the whole box $w(A)\times A$ lies in $\lambda_C$. Hence, for every $\sigma\in\bS_k$, replacing the points of $T(w)$ in the columns $A$ by those of $w_{A\leftarrow\sigma}$ gives a transversal of $\lambda_C$, and by Proposition~\ref{prop:classes}(1) we have $w_{A\leftarrow\sigma}\in C$. So all terms of $e_k[w,A]$ lie in $C$, and under the isomorphism \eqref{eq:class-dec} the element $e_k[w,A]$ becomes $e_k[T(w),A]$, one of the elements spanning $(e_k)_{\lambda_C}$.

Second, let $C$ be a $\cW$-class, $T$ a transversal of $\lambda_C$, and $A$ a set of $k$ columns such that the corner $c$ of the points of $T$ in the columns $A$ lies in $\lambda_C$. Let $w\in C$ be the permutation with $T(w)=T$. Then $c\in\lambda_\cW(w)$, that is, $c$ has a $\cW$-witness in $w$. So $e_k[w,A]$ is one of the elements of Lemma~\ref{lem:bwx-ideal}, and it corresponds to $e_k[T,A]$.

Finally, by the first part every spanning element of $\cI_{k,\cW}(N)$ lies in a single summand $\kk C$. So $\cI_{k,\cW}(N)$ is the direct sum of the subspaces $\cI_C$, each spanned by the spanning elements lying in $\kk C$, and by the first two parts $\cI_C$ corresponds to $(e_k)_{\lambda_C}$, which is $I_k(\lambda_C)_{(1;1)}$ by Lemma~\ref{lem:minors}. Taking quotients gives \eqref{eq:splitting}.
\end{proof}

\begin{figure}[ht]
\centering
\begin{tikzpicture}[scale=0.55]
\shaderows{1/5,2/5,3/5,4/1,5/1}
\gridlines{6}{6}
\pt{5}{1}\pt{1}{3}\pt{2}{4}\pt{3}{5}
\opt{6}{2}\opt{4}{6}
\end{tikzpicture}
\caption{The $\cW$-board of $u=561234$ for $\cW=\{1\}$. The $\cW$-inner points are black and the $\cW$-outer points are white.}\label{fig:tau1}
\end{figure}
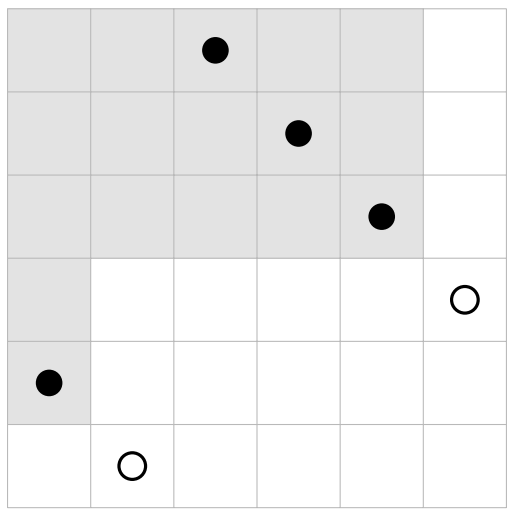

\begin{example}\label{ex:splitting}
Let $\cW=\{1\}$ and $u=561234$ (Figure~\ref{fig:tau1}). The right-to-left maxima are $(6,2)$ and $(4,6)$, and $\lambda_\cW(u)$ has row lengths $5,5,5,1,1,0$. The $\cW$-inner points are $(5,1)$, $(1,3)$, $(2,4)$ and $(3,5)$, and the board $\lambda_C\subset\{1,2,3,5\}\times\{1,3,4,5\}$ has row lengths $4,4,4,1$. Its transversals send column $1$ to row $5$ and the columns $3,4,5$ to the rows $1,2,3$ in any order, so the class of $u$ is $C=\{56abc4:\{a,b,c\}=\{1,2,3\}\}$, of size six. For $k=2$ and $A=\{3,4\}$ the corner $(2,4)$ has the $\cW$-witness $\{(4,6)\}$, and
\[
e_2[u,\{3,4\}]=(e_2\oplus1)[u,\{3,4,6\}]=561234-562134=\Delta_{\{1,2\},\{3,4\}}\,x_{51}x_{62}x_{35}x_{46}.
\]
Under the isomorphism $\kk C\cong\kk\lambda_C$ this element becomes the minor of $\lambda_C$ on the rows $1,2$ and the columns $3,4$, times $x_{51}x_{35}$.
\end{example}

\begin{example}\label{ex:tau1count}
Let $\cW=\{1\}$ and $N=4$. By Example~\ref{ex:classes} there are $14$ classes: eight of size one, five of size two, and the class $\{abc4\}$ of size six, whose board is the full $3\times3$ square.
\begin{itemize}
\item For $k=2$, every class contributes one dimension to \eqref{eq:splitting} by Remark~\ref{rem:k2board}. This gives $\dim(\cS/(e_2\oplus1))(4)=14=\#\Av_4(213)=\#\Av_4(123)$.
\item For $k=3$, the classes of size at most two contain no $3$-chains and contribute their size. The class of size six contributes $6-1=5$, since the $3\times3$ determinant is its only relation. This gives $\dim(\cS/(e_3\oplus1))(4)=8+10+5=23=\#\Av_4(3214)=\#\Av_4(1234)$.
\end{itemize}
The class $\{3412\}$, whose board is a staircase, illustrates Remark~\ref{rem:projection}. For $k=2$ the ladder ideal of this board vanishes, and the class contributes one dimension; indeed $3412$ avoids both $213$ and $123$. The intersection of the rank variety with the matrices vanishing outside the board would kill this class.
\end{example}

\begin{remark}\label{rem:strata}
	The splitting \eqref{eq:splitting} is the shadow of a stratification. Let $J_N\subset R_N$ be the ideal generated by the products $\Delta_{A,B}\,x_X$, where $\Delta_{A,B}$ is a $k\times k$ minor and $X$ is a set of cells strictly south-east of its corner, in distinct rows and columns, forming a pattern from $\cW$. By Lemmas~\ref{lem:minors} and~\ref{lem:bwx-ideal}, its weight space $(J_N)_{(1;1)}$ is $\cI_{k,\cW}(N)$. For a matrix $M$ let $\lambda_\cW(M)$ be the set of cells $c$ such that the nonzero entries of $M$ contain such a set $X$ strictly south-east of $c$. Then $M$ lies in the zero set of $J_N$ if and only if $\operatorname{rank}M_{[a]\times[b]}<k$ for all $(a,b)\in\lambda_\cW(M)$, that is, if and only if the restriction of $M$ to its own board $\lambda_\cW(M)$ lies in the ladder determinantal variety $Y_k(\lambda_\cW(M))$. The proofs of Lemma~\ref{lem:bw} and Proposition~\ref{prop:classes}, applied to the supports of monomials instead of permutations, show that $\lambda_\cW(M)$ only depends on the set $O$ of nonzero entries of $M$ outside $\lambda_\cW(M)$, and that there is a decomposition of multigraded vector spaces
	\[
	R_N/J_N\cong\bigoplus_O\ x_O\,\kk[x_c:c\in O]\otimes R_{\lambda_O}/I_k(\lambda_O),
	\]
	where $O$ runs over the sets of cells of $[N]\times[N]$ disjoint from their own board $\lambda_O$. Geometrically, the zero set of $J_N$ is the disjoint union of the strata $(\kk^\times)^O\times Y_k(\lambda_O)$. In multidegree $(1;1)$, the set $O$ becomes a partial permutation, the set of $\cW$-outer points of a class $C$, and the summand becomes $(R_{\lambda_C}/I_k(\lambda_C))_{(1;1)}$.
\end{remark}

\subsection{The Gr\"obner bases}\label{sec:BWX-proof}

\begin{lemma}\label{lem:bwx-chain}
Let $u\in\bS_N$, let $C$ be its $\cW$-class, and let $p$ be $k\cdots21$ or $12\cdots k$. Then $u$ contains a pattern from $p\oplus\cW$ if and only if the transversal $T(u)$ contains $p$ in $\lambda_C$.
\end{lemma}

\begin{proof}
An occurrence of $p\oplus\tau$ in $u$ consists of an occurrence of $p$ and of an occurrence of $\tau$ all of whose points lie below and to the right of all points of the former, that is, strictly south-east of its corner. So $u$ contains a pattern from $p\oplus\cW$ if and only if $u$ has an occurrence of $p$ whose corner has a $\cW$-witness in $u$, that is, lies in $\lambda_\cW(u)$. By Lemma~\ref{lem:bw}(3), the points of such an occurrence are $\cW$-inner, so they belong to $T(u)$, and their corner lies in $\rho_C\times\kappa_C$, hence in $\lambda_C$. Conversely, if some points of $T(u)$ form the pattern $p$ and their corner lies in $\lambda_C\subset\lambda_\cW(u)$, then this corner has a $\cW$-witness in $u$.
\end{proof}

\begin{proof}[Proof of Theorem~\ref{thm:bwx}]
First, consider $<_\lex$. By Theorem~\ref{thm:splitting}, every element of $\cI_{k,\cW}(N)$ is a sum of elements of the spaces $\cI_C$, and these spaces are supported on disjoint sets of permutations. So the leading permutations of $\cI_{k,\cW}(N)$ are the leading permutations of all the $\cI_C$. By Theorem~\ref{thm:splitting} and Proposition~\ref{prop:classes}(2), the leading permutations of $\cI_C$ correspond to the leading transversals of $(e_k)_{\lambda_C}$, and by Corollary~\ref{cor:shape-wilf} these are the transversals of $\lambda_C$ containing $k\cdots21$ in $\lambda_C$. By Lemma~\ref{lem:bwx-chain}, the corresponding permutations are the elements of $C$ containing a pattern from $k\cdots21\oplus\cW$. Hence every leading permutation of $\cI_{k,\cW}$ contains the leading monomial of one of the generators, that is, the generators form a Gr\"obner basis for $<_\lex$.

Second, the same argument with $12\cdots k$ in place of $k\cdots21$ applies to $<_\lex^c$.

Finally, the sign twist $\theta$ maps $e_k\oplus\tau$ to $\sgn(\tau)\,s_k\oplus\tau$ and preserves leading monomials, so $s_k\oplus\cW$ is a Gr\"obner basis for the same orderings. The equality of the numbers of avoiders is Proposition~\ref{prop:wilf}.
\end{proof}

The proof shows more: for every $\cW$-class $C$, the image of $\kk C$ in the quotient has a basis of the transversals of $\lambda_C$ without a decreasing $k$-chain in $\lambda_C$, and another basis of those without an increasing $k$-chain in $\lambda_C$. Summing over $\Cl_\cW(N)$ gives the reduction of Backelin, West and Xin of their theorem to the shape-Wilf equivalence of $12\cdots k$ and $k\cdots21$.

\begin{remark}\label{rem:board-general}
The only property of $e_k$ used in Sections~\ref{sec:splitting} and~\ref{sec:BWX-proof} is that it is a board Gr\"obner basis. Let $G$ be any board Gr\"obner basis for $<_\lex$ or for $<_\lex^c$ (Definition~\ref{def:board-gb}). Since $(g\oplus\tau)[w,L]=g[w_{L\leftarrow\id\oplus\tau},A]$, where $A$ consists of the first $\operatorname{arity}(g)$ elements of $L$, the same proofs show that $\{g\oplus\tau:g\in G,\ \tau\in\cW\}$ is a Gr\"obner basis for the same ordering, and that the quotient by the ideal it generates is $\bigoplus_{C\in\Cl_\cW(N)}\kk\lambda_C/(G)_{\lambda_C}$ in arity $N$.
\end{remark}

\begin{remark}\label{rem:binomial-fails}
For $k\ge3$ the sign-twisted sum $e_k\oplus\tau$ cannot be replaced by the binomial $12\cdots k\oplus\tau-k\cdots21\oplus\tau$. For instance, by Proposition~\ref{prop:binomial}(1) the quotient by the ideal generated by $1234-3214$ has dimension $99$ in arity five, while $\#\Av_5(1234)=\#\Av_5(3214)=103$.
\end{remark}

\end{document}